\documentstyle[12pt]{article}

\input amssym.def
\input amssym
\begin{document}
\def \Z{\Bbb Z}
\def \C{\Bbb C}
\def \N{\Bbb N}
\def \1{{\bf 1}}
\def \L{{\cal L }}
\def \U{{\cal U }}
\def \V{{\cal V }}
\def \e{\varepsilon }
\def \Res{\mbox{Res}}

\def \pf{\noindent {\bf Proof.\ }}
\def \qed{\ \ \ $\Box$}

\newtheorem{thm}{Theorem}[section]
\newtheorem{prop}[thm]{Proposition}
\newtheorem{coro}[thm]{Corollary}
\newtheorem{conj}[thm]{Conjecture}
\newtheorem{lem}[thm]{Lemma}
\newtheorem{rem}[thm]{Remark}
\newtheorem{de}[thm]{Definition}
\makeatletter
\@addtoreset{equation}{section}
\def\theequation{\thesection.\arabic{equation}}
\makeatother
\makeatletter

\begin{center}{\Large \bf Vertex algebras generated by Lie algebras}

\vspace{0.5cm}
Mirko Primc
\\ Department of Mathematics, University of Zagreb, Croatia
\end{center}

\begin{abstract}
In this paper we introduce a notion of vertex Lie algebra $U$, in a way a
``half''
of vertex algebra structure sufficient to construct the corresponding
local Lie algebra $\L(U)$ and a vertex algebra $\V(U)$. We show that we may
consider
$U$ as a subset $U \subset {\cal V}(U)$ which generates ${\cal V}(U)$ and
that
the vertex Lie algebra structure on $U$ is induced by the vertex algebra
structure on
$\V(U)$. Moreover, for any vertex algebra $V$ a given homomorphism $U
\rightarrow V$
of vertex Lie algebras extends uniquely to a homomorphism ${\cal V}(U)
\rightarrow V$
of vertex algebras. In the second part of paper we study under what
conditions on
structure constants one can construct a vertex Lie algebra $U$ by starting
with a given
commutator formula for fields.
\end{abstract}

\addtocounter{footnote}{1}%
\footnotetext{Supported in part by the Ministry of Science and Technology of
the Republic of Croatia, grant 037002.}
\addtocounter{footnote}{1}%
\footnotetext{1991 Mathematics Subject Classification: 17B69.}

\section{Introduction}

One way of seeing a vertex algebra $V$ is as a vector space with infinitely
many bilinear multiplications $u_n v$, $n\in\Z$, which correspond to
``normal order products'' $u(z)_n v(z)$ of fields
$u(z)$ and $v(z)$ associated to vectors $u$ and $v$. For two fields there is
a formula in which the commutator is expressed in terms of
products  $u(z)_n v(z)$ only for $n\geq 0$. So any vector space $U\subset V$
closed for
``positive'' multiplications will look like some kind of Lie algebra.

In this paper we define a vertex Lie algebra as a vector space $U$ given
infinitely
many bilinear multiplications $u_n v$, $n\in\N$, satisfying a Jacobi
identity
and a skew symmetry in terms of given derivation $D$.
In particular, any vertex algebra is a vertex Lie algebra if we forget the
multiplications other than for $n\geq 0$. This structure is
sufficient to construct a Lie algebra ${\cal L}(U)$ generated by vectors
$u_n$\,, $u\in U$, $n\in \Z$, such that the commutator formula for fields
with $u(z)=\sum_{n\in \Z} u_n z^{-n-1}$ defines the commutator in ${\cal
L}(U)$.
This Lie algebra has the obvious decomposition
${\cal L}(U)={\cal L}_-(U)\oplus {\cal L}_+(U)$ and the induction by a
trivial ${\cal L}_+(U)$-module gives a generalized Verma ${\cal
L}(U)$-module
\[
   {\cal V}(U)={\cal U}\left({\cal L}(U)\right)
\otimes_{\,{\cal U}\left({\cal L}_+(U)\right)} \C \cong {\cal U}\left({\cal
L}_-(U)\right),
\]
where $\cal U$ stands for the universal enveloping algebra of a given Lie
algebra.

We show that ${\cal V}(U)$ is a vertex algebra and that we may consider
$U$ as a subset $U \subset {\cal V}(U)$ which generates ${\cal V}(U)$.
So an ``abstract'' $U$ turns to be as in the
``concrete'' motivating example at the beginning. Moreover,
for any vertex algebra $V$ a given homomorphism $U \rightarrow V$ of vertex
Lie
algebras extends uniquely to a homomorphism ${\cal V}(U) \rightarrow V$ of
vertex
algebras. Because of this universal property we call ${\cal V}(U)$ the
universal enveloping vertex algebra of $U$. These constructions are modeled
after
and apply to the well known examples of vertex (super)algebras associated to
affine Lie
algebras, Virasoro algebra and Neveu-Schwarz algebra.

In the second part of this paper we study under what conditions on structure
constants one
can construct a vertex algebra by starting with a given commutator formula
for fields, or equivalently, by starting with a given singular part of
operator product expansion for fields.
In the case of a commutator formula closed for a set $S$ of quasi-primary
fields we give explicit necessary and sufficient conditions for the
existence of
universal vertex algebra $\V(\langle S\rangle)$ generated by $S$.

This work rests on the inspiring results of Hai-sheng Li in [Li] and in a
way
a complementary theorem on generating fields in [FKRW] and [MP], but it goes
without saying
that many ideas used here stem from [B], [FLM], [FHL], [DL] and [G], to
mention just
a few.
Li's point of view on generating fields gives a natural framework for
studying
modules and it was tempting to see how some of his Lie-theoretic arguments
could be extended to a more general setting.
The key technical point
is the observation that a direct proof that the local algebra ${\cal L} (V)$
of
vertex algebra $V$ is a Lie algebra involves just a ``positive half'' of
both
the Jacobi identity and the skew symmetry, i.e., to be more precise,
involves
only the principal part of the formal Laurent series which appear in
these relations. In a similar way only a ``positive half'' of
the commutator formula implies the ``positive half'' of the Jacobi identity.
Some of these arguments are essentially
a simpler copy of Li's arguments for vertex (super)algebras and work just
the same
in the vertex superalgebra case.

After this work was completed I have learned that the notion of vertex Lie
algebra is a special case of a more general notion of local vertex Lie
algebra introduced and studied in [DLM] by C. Dong, H.-S. Li and G. Mason.
In part their work contains some results similar to the ones described
above.
I thank the authors for sending me their paper and for informing me about
some references in physics literature.

I thank Arne Meurman, Jim Lepowsky and Haisheng Li for different
ways in which they greatly contributed to this work.

\section{Generating fields for vertex algebras}

For a $\Z _2$-graded vector space $W=W^0+W^1$ we write $|u|\in\Z_2$, a
degree of
$u$, only for homogeneous
elements: $|u|=0$ for an even element $u\in W^0$ and $|u|=1$ for an odd
element
$u\in W^1$. For any two $\Z _2$-homogeneous elements $u$ and $v$ we define
$\e _{u,v}=(-1)^{|u||v|}\in \Z$.

By following [Li], a vertex superalgebra $V$ is a $\Z _2$-graded vector
space
$V=V^0+V^1$ equipped with a specified vector \1 called the vacuum vector, a
linear
operator $D$ on $V$ called the derivation and a linear map
$$
V\rightarrow (\mbox{End}\,V)[[z^{-1},z]],\ \ \ \
v\mapsto Y(v,z)=\sum_{n\in\Z}v_nz^{-n-1}
$$
satisfying the following conditions for $u,v\in V$:
\begin{eqnarray}
& &u_nv=0\ \ \
\mbox{for}\ \ n \ \ \mbox{sufficiently\ large};\label{vlw0}\\
& &[D,Y(u,z)]=Y(Du,z)=\frac{d}{dz}Y(u,z);\label{der}\\
& &Y({\bf 1},z)=\mbox{id}_V\ \ \ (\mbox{the identity operator on}\
V);\label{vacuum}\\
& &Y(u,z)\1 \in (\mbox{End}\,V)[[z]]\ \ \ \mbox{and}\ \ \ \lim_{z\to
0}Y(u,z)\1 =u;
\label{creation}
\end{eqnarray}
For $\Z _2$-homogeneous elements $u,v\in V$ the Jacobi identity holds:
 \begin{equation}\label{jacobi}
\begin{array}{c}
\displaystyle{z^{-1}_0\delta\left(\frac{z_1-z_2}{z_0}\right)
Y(u,z_1)Y(v,z_2)-\e _{u,v}\, z^{-1}_0\delta\left(\frac{z_2-z_1}{-z_0}\right)
Y(v,z_2)Y(u,z_1)}\\
\displaystyle{=z_2^{-1}\delta\left(\frac{z_1-z_0}{z_2}\right)
Y(Y(u,z_0)v,z_2)}.
\end{array}
\end{equation}
Finally, for any $\Z _2$-homogeneous $u,v\in V$ and $n\in \Z$ we assume
\begin{equation}\label{homog}
    |u_nv|=|u|+|v|
\end{equation}
(i.e. $u_n v$ is homogeneous and $|u_nv|=|u|+|v|$).

Sometimes we will emphasize that $Y(u,z)$ is pertinent to a vertex
superalgebra $V$ by writing $Y_V(u,z)$.

In the definition of vertex superalgebra, VSA for short, the condition
(\ref{creation}) can be equivalently replaced by the condition
$ Y(u,z)\1 =e^{zD}u$ for all $u\in V$, and both are called the creation
property. In a way the creation property is a special case of the
skew symmetry
$$
Y(u,z)v=\e_{u,v}e^{zD}Y(v,-z)u
$$
which holds for all homogeneous $u,v\in V$.

Note that in the case when $V=V^0$,
i.e. when all vectors are even, all $\Z_2$-grading conditions become trivial
and we speak of a vertex algebra, VA for short.
In general, for a $\Z _2$-graded vector space $W$ the vector spaces
$$\begin{array}{l}
\mbox{End}\,W=(\mbox{End}\,W)^0\oplus (\mbox{End}\,W)^1,\\
(\mbox{End}\,W)[[z^{-1},z]]=(\mbox{End}\,W)^0 [[z^{-1},z]]
\oplus (\mbox{End}\,W)^1[[z^{-1},z]]
\end{array}$$
are $\Z_ 2$-graded as well and our assumption (\ref{homog}) implies
that for a homogeneous element $u\in V$ the operators $u_n\,$ are
homogeneous
for all $n\in \Z$, that $Y(u,z)$ is homogeneous and
$$
|u|=|u_n|=|Y(u,z)|.
$$
This together with (\ref{creation}), (\ref{vacuum}) and
(\ref{der}) implies that $\1 \in V$ is even and that $D\in \mbox{End}\,V$ is
even.
On the other hand, $\1 \in V$ is even together with
(\ref{vlw0}), (\ref{creation}) and (\ref{jacobi}) implies (\ref{der}) for
$D$ defined by $Du=u_{-2}\1$; the relation (\ref{vacuum}) follows as well.
For this reason we define a homomorphism $\varphi : V\rightarrow U$ of
two vertex superalgebras as a $\Z_2$-grading preserving linear map $\varphi$
such that
$$
\varphi (u_n v)=(\varphi (u))_n(\varphi (v))
$$
for all $u,v\in V$, $n\in \Z$; as a consequence we have relations
$$
\varphi (\1)= \1,\ \ \ \varphi D=D\varphi.
$$

For a subset $U\subset V$ we denote by $\langle U \rangle$ a vertex
superalgebra
generated by $U$, i.e. the smallest vertex superalgebra containing the set
$U$.

A module $M$ for a vertex superalgebra $V$ is a $\Z _2$-graded vector space
$M=M^0+M^1$ equipped with an even linear operator $D\in (\mbox{End}\,M)^0$
and a linear map
$$
V\rightarrow (\mbox{End}\,M)[[z^{-1},z]],\ \ \ \
v\mapsto Y_M(v,z)=\sum_{n\in\Z}v_nz^{-n-1}
$$
satisfying the following conditions for $u,v\in V$ and $w\in M$:
\begin{eqnarray}
& &u_nw=0\ \ \
\mbox{for}\ \ n\in \Z \ \ \mbox{sufficiently large};\label{vlw0M}\\
& &[D,Y_M(u,z)]=Y_M(Du,z)=\frac{d}{dz}Y_M(u,z);\label{derM}\\
& &Y_M({\bf 1},z)=\mbox{id}_M\ \ \ (\mbox{the identity operator on}\
M);\label{vacuumM}
\end{eqnarray}
For $\Z _2$-homogeneous elements $u,v\in V$ the Jacobi identity holds:
 \begin{equation}\label{jacobiM}
\begin{array}{c}
\displaystyle{z^{-1}_0\delta\left(\frac{z_1-z_2}{z_0}\right)
Y_M(u,z_1)Y_M(v,z_2)-\e _{u,v}\,
z^{-1}_0\delta\left(\frac{z_2-z_1}{-z_0}\right)
Y_M(v,z_2)Y_M(u,z_1)}\\
\displaystyle{=z_2^{-1}\delta\left(\frac{z_1-z_0}{z_2}\right)
Y_M(Y(u,z_0)v,z_2)}.
\end{array}
\end{equation}
Finally, for any $\Z _2$-homogeneous $u\in V$, $w\in M$ and $n\in \Z$ we
assume
\begin{equation}\label{homogM}
    |u_nw|=|u|+|w|
\end{equation}
(i.e. $u_n w$ is homogeneous and $|u_nw|=|u|+|w|$).

Clearly, $V$ is a $V$-module with $Y_V(u,z)=Y(u,z)$, and, as before, for a
$V$-module
$M$ we have
$$
|u|=|u_n|=|Y_M(u,z)|.
$$
Let $M$ be a $V$-module and $u,v\in V$ homogeneous elements. Set
$u(z)=Y_M(u,z)$ and $v(z)=Y_M(v,z)$. As a consequence of the definition of
$V$-module $M$ we have the normal order product formula and the locality:
The normal order product formula states that the field $Y_M(u_nv,z)$ equals
\begin{equation}\label{product}
   u(z)_nv(z)=\Res_{z_1}\left((z_1-z)^n u(z_1)v(z)-
(-1)^{|u(z)||v(z)|}(-z+z_1)^n v(z)u(z_1)\right).
\end{equation}
The locality property states that for some $N=N(u,v)\in \N$
\begin{equation}\label{locality}
   (z_1-z_2)^N u(z_1)v(z_2)=(-1)^{|u(z)||v(z)|}(z_1-z_2)^N v(z_2)u(z_1).
\end{equation}
In the case when formal Laurent series $u(z)$ and $v(z)$ satisfy
relation (\ref{locality}) for some $N=N(u(z),v(z))\in \N$, we say
that $u(z)$ and $v(z)$ are local to each other.

Hai-sheng Li proved [Li, Proposition 2.2.4] that in the definition of
vertex superalgebra the Jacobi identity  (\ref{jacobi})  can be
equivalently substituted by the locality property (\ref{locality}).

Let $M$ be a $\Z _2$-graded vector space
$M=M^0+M^1$ equipped with an even linear operator $D\in (\mbox{End}\,M)^0$.
By following [Li] define a vertex operator on $M$ as
a homogeneous formal Laurent series $v(z)=\sum_{n\in\Z}v_nz^{-n-1}$ in
$(\mbox{End}\,M)[[z^{-1},z]]$ such that the property (\ref{vlw0M}) holds
for all $w\in W$, that $v(z)$ is local with itself (i.e. that
(\ref{locality}) holds for $u(z)=v(z)$) and that
\begin{equation}\label{derv(z)}
  [D,v(z)]=\frac{d}{dz}\,v(z).
\end{equation}
We shall also say that a linear combination of vertex operators on $M$ is a
vertex operator on $M$.

It is clear that giving a map $Y(\cdot,z)$ is equivalent to giving
infinitely
many bilinear multiplications $u_nv$, $n\in \Z$. Let $F(M)$ be a space of
formal
Laurent series $v(z)$ in $(\mbox{End}\,M)[[z^{-1},z]]$
such that the property (\ref{vlw0M}) holds for all $w\in W$. Then $F(M)$ is
$\Z_2$-graded and on $F(M)$ bilinear multiplications $ u(z)_nv(z)$, $n\in
\Z$,
given by (\ref{product}) are well
defined. Moreover, Hai-sheng Li proved the following theorem
 [Li, Corollary 3.2.11 of Theorem 3.2.10]:

\begin{thm}\label{LiThm}
Let $M$ be any $\Z_2$-graded vector space
equipped with an even linear operator $D$ and let $U$ be any set of
mutually local homogeneous vertex operators on $M$. Let $\langle U\rangle$
be
the subspace of $F(M)$ generated by $U$ and $I(z)=\mbox{id}_M$ under the
vertex
operator multiplication (\ref{product}). Then $\langle U\rangle$ is a vertex
superalgebra with the vacuum vector \ $\1=I(z)$ and the derivation
$D=\frac{d}{dz}$.
Moreover, $M$ is a $\langle U\rangle$-module.
\end{thm}

Clearly this theorem implies that for a vertex superalgebra $V$ the set
of fields $\{Y(u,z)\mid u\in V\}$ is a vertex superalgebra, and by
construction
it is clear that $u\mapsto Y(u,z)$ is an isomorphism. Because of this
isomorphism we
shall sometimes say that for a subset $U\subset V$ the vertex superalgebra
$\langle U \rangle\subset V$ is generated by the set of fields
$\{Y(u,z)\mid u\in U\}$.

For the following theorem
cf.\ [G], [Xu], [FKRW] and [MP2]. The theorem was
proved in [MP] for graded vertex
algebras, the results in [Li] allow us to extend the proof to the vertex
superalgebra case:
\begin{thm}\label{MPThm}
Let V be a $\Z _2$-graded vector space
$V=V^0+V^1$ equipped with an even linear operator $D\in (\mbox{End}\,V)^0$
and
an even vector $\1$ such that $D\1=0$.
Let $U$ be a $\Z_2$-graded subspace of $V$ given a linear map
$$
Y : U \rightarrow (\mbox{End}\,V)[[z,z^{-1}]],\ \ \  u \mapsto Y(u,z) =
\sum_{n \in \Z}
u_n z^{-n-1}
$$
such that $\{Y(u,z)\mid u\in U^0\cup U^1\}$ is a set of mutually local
vertex
operators on $V$ satisfying the following two conditions:
\begin{equation}\label{2.32}
Y(u,z)\1 \in (\mbox{End}\,V)[[z]]\ \ \ \mbox{and}\ \ \ \lim_{z\to 0}Y(u,z)\1
=u;
\end{equation}
\begin{equation}\label{2.34}
V =  \mbox{span\,} \{u^{(1)}_{n_1} \cdots u^{(k)}_{n_k} \1 \mid  k \in \N ,
n_i \in \Z, u^{(i)} \in U^0\cup U^1  \} .
\end{equation}
Then  $Y$ extends uniquely into a vertex superalgebra with the vacuum vector
$\1$
and the derivation $D$.
\end{thm}

\pf
Let $W$ be the $\Z_2$-graded space of all vertex operators
$a(z)$ on $V$ such that
$a(z)$ and $Y(u,z)$ are mutually local for each $u\in U^0\cup U^1$
and that $a(z)\1$ is a power series in $z$. Define a linear map
$$
{\phi}: W\rightarrow V, \ \ \ a(z)\mapsto a_{-1}\1.
$$
Note that our assumptions imply $Y(u,z)\in W$ for all $u\in U$.
Hence  by (\ref{2.32})
$$
U\subset {\phi}(W)\subset V .
$$
Also note that $\phi$ preserves the $\Z_2$-grading and that $\phi (W)\subset
V$ is
a $\Z_2$-graded subspace.

\noindent{\bf Step 1.}\ \ {\it ${\phi}$ is injective.}
Assume that $a(z)\in W$ is homogeneous and $a_{-1}\1 = 0$. Since $a(z)\1$ is
a power series
in $z$, (\ref{derv(z)}) and $D\1=0$ imply
\begin{equation}\label{2.36}
a(z)\1 = e^{zD}a_{-1}\1 = 0.
\end{equation}
Let $X = \{v\in V\mid a(z)v = 0\}$ and let $v\in X$, $u\in U^0\cup U^1$.
Then there
exists $N\in \N$ such that
$$
(z_1-z_2)^Na(z_1)Y(u,z_2)v = (-1)^{|a(z)||u|}(z_1-z_2)^NY(u,z_2)a(z_1)v= 0.
$$
Hence $a(z_1)u_nv = 0 $ and $u_nX\subset X $. Since
(\ref{2.36}) implies $\1\in X $, it follows from (\ref{2.34})
that $X = V $. Hence $a(z) = 0 $ and ${\phi}$ is injective.

\noindent{\bf Step 2.}  Set $X = {\phi}(W)$. Since ${\phi}$ is injective we
can define
$$
Y: X\to(\mbox{End}\, V)[[z,z^{-1}]]\ \ \ \mbox{by}\ \ \ Y(v,z) =
{\phi}^{-1}(v)
$$
for $v\in X $. By (\ref{2.32}) the two meanings of $Y(u,z)$ for $u\in U$
denote the
same series. {\it We claim that $X = {\phi}(W) = V$,
i.e. we have well defined vertex operators $Y(v,z) $ for all $v\in V $.}

For $n\in\Z$ and homogeneous elements $u\in U $, $v\in X $, define
$Y(u_nv,z) $ by
using the product (\ref{product}):
$$
Y(u_nv,z) = \left(Y(u,z)\right)_n\left(Y(v,z)\right).
$$
Since $Y(u,z)$ and $Y(v,z)$ are in $W$, it follows from
[Li, Lemma 3.1.4, Lemma 3.1.8 and Proposition 3.2.7] that $Y(u_nv,z)$ is a
vertex
operator on $V$. Again by [Li, Proposition 3.2.7]
vertex operators $Y(u_nv,z)$ and $Y(w,z)$ are
mutually local for each $w\in U^0\cup U^1$. It follows from (\ref{product})
that $Y(u_nv,z)\1$ is a power series in $z$. Hence  $Y(u_nv,z)\in W $.
It follows from (\ref{product}) (cf. [Li, (3.1.9)]) that
$$
{\phi}\left(Y(u_nv,z)\right)=(u_nv)_{-1}\1 = u_nv_{-1}\1 = u_nv.
$$
Hence $u\in U $, $v\in X $, $n\in \Z $ implies
\begin{equation}\label{2.37}
u_nv\in X.
\end{equation}
Since $\mbox{id}_Vz^0\in W $ and $\mbox{coeff}_{z^0}(\mbox{id}_Vz^0)\1 = \1
$, we have
\begin{equation}\label{2.38}
\1\in X
\end{equation}
and $Y(\1,z) = \mbox{id}_Vz^0=\mbox{id}_V$.
Now (\ref{2.34}), (\ref{2.37}), (\ref{2.38}) imply $X = V$,
i.e. we have well defined vertex operators $Y(v,z) $ for all $v\in V $.

\noindent {\bf Step 3.}\ \ {\it $Y(u,z)$ and $Y(v,z) $ are mutually local
for all pairs
$u,v\in V $.}
Fix a homogeneous element $v^{(0)}\in V $ and set
$$X =\mbox{span}\,\{v\in V\mid Y(v,z) \mbox{ and }Y(v^{(0)},z)
\mbox{ are mutually local}\}.
$$
Clearly $\{\1\}\cup U\subset X$. For $n\in\Z $ and homogeneous elements
$u\in U $, $v\in X $ by definition $Y(u_nv,z) =
\left(Y(u,z)\right)_n\left(Y(v,z)\right)$,
so [Li, Proposition 3.2.7] implies that  $u_nv\in X $. Hence (\ref{2.34})
again
implies $X = V$ as required.

We have showed that $V$ satisfies the conditions of [Li, Proposition 2.2.4],
so
$V$ is a vertex superalgebra. \qed


\section{Vertex Lie algebras}

Let $V$ be a vertex superalgebra. By taking
$\mbox{Res}_{z_0}\mbox{Res}_{z_1}\,z_0^n $ of the
Jacobi identity for $V$ or the Jacobi identity for a $V$-module $M$,
that is the coefficients of $z_0^{-n-1}z_1^{-1}$, we get the
normal order product formula (\ref{product}). These relations
are the components of associator formula obtained by taking the
residue $\mbox{Res}_{z_1}$ of the Jacobi identity:
\begin{eqnarray}
&& Y(Y(u,z_0)v,z_2)  - Y(u,z_0+z_2)Y(v,z_2)\nonumber\\
&=&-\e _{u,v}\,\mbox{Res}_{z_1}
z^{-1}_0\delta\left(\frac{z_2-z_1}{-z_0}\right)
Y(v,z_2)Y(u,z_1)\nonumber\\
&=&\e_{u,v}\,Y(v,z_2)\{Y(u,z_2+z_0) - Y(u,z_0+z_2)\}.\nonumber
\end{eqnarray}
By taking the residue $\mbox{Res}_{z_0} $ of the Jacobi identity we get
the commutator formula
\begin{eqnarray}
 [Y(u,z_1), Y(v,z_2)]  &= & \mbox{Res}_{z_0} z_2^{-1} \delta
\left(\frac{z_1 - z_0}
{z_2} \right) Y(Y(u,z_0)v,z_2) \label{commutator}\\
        &= & \sum_{i \geq 0} \frac{(-1)^i}{i!} \left(\frac{d}{dz_1}
\right)^i z^{-1}_2
\delta(z_1/z_2)Y(u_iv,z_2) .\nonumber
\end{eqnarray}
We can write these identities for components $u_m$ and $v_n$ of $Y(u,z)$ and
$Y(v,z)$:
$$
u_mv_n-\e_{u,v}\,v_nu_m = \sum_{i \geq 0} {m \choose i}(u_iv)_{n+m-i} \ .
$$

In general, if we take
$\mbox{Res}_{z_0}\mbox{Res}_{z_1}\mbox{Res}_{z_2}\,z_0^k z_1^m z_2^n$
of the Jacobi identity applied to a vector $w$, that is the coefficients of
$z_0^{-k-1}z_1^{-m-1}z_2^{-n-1}$, we get for components of vertex
operators the identities
\begin{equation}\label{jacobicomp}
\begin{array}{c}
\displaystyle{\sum_{i\geq0} (-1)^i {k \choose i}(u_{m+k-i}(v_{n+i}w) -
\e_{u,v}\,(-1)^k v_{n+k-i}(u_{m+i}w)) }\\
\displaystyle{=\sum_{i \geq 0} {m \choose i}(u_{k+i}v)_{m+n-i}w}.
\end{array}
\end{equation}
These relations hold for all $k, m, n\in \Z$, but it should be noticed that
for $k, m, n\in \N$ these relations involve only indices in $\N$. In a way
the purpose of this paper is to study the consequences of this ``half'' of
the Jacobi identity.

Let $A$ and $B$ be two formal Laurent series (in possibly several variables
$z_0, z_1,\dots$). We shall write $A\simeq B$ if the principal parts of $A$
and $B$
are equal. For example, $A(z_0,z_1,z_2)\simeq B(z_0,z_1,z_2)$ means that the
coefficients of $z_0^{-k-1}z_1^{-m-1}z_2^{-n-1}$ in $A$ equal the
coefficients in $B$ for all $k, m, n\in \N$. In particular,
the set of relations (\ref{jacobicomp}) for all $k, m, n\in \N$ is
equivalent to
the ``half Jacobi identity'' (\ref{Hjacobi}) written below. In a similar way
we shall
speak of the half commutator
formula or the half associator formula. From the way they were obtained, it
is clear they are a subset of the half Jacobi identity viewed by components.

With the above notation we define a vertex Lie superalgebra $U$ as a $\Z
_2$-graded
vector space $U=U^0+U^1$ equipped with an even linear
operator $D$ on $U$ called the derivation and a linear map
$$
U\rightarrow z^{-1}(\mbox{End}\,U)[[z^{-1}]],\ \ \ \
v\mapsto Y(v,z)=\sum_{n\geq 0}v_nz^{-n-1}
$$
satisfying the following conditions for $u,v\in U$:
\begin{eqnarray}
& &u_nv=0\ \ \
\mbox{for}\ \ n \ \ \mbox{sufficiently\ large};\label{Hvlw0}\\
& &[D,Y(u,z)]=Y(Du,z)=\frac{d}{dz}Y(u,z);\label{Hder}
\end{eqnarray}
For $\Z _2$-homogeneous elements $u,v\in V$ the half skew symmetry holds:
\begin{equation}\label{HSS}
Y(u,z)v\simeq \e_{u,v}e^{zD}Y(v,-z)u.
\end{equation}
For $\Z _2$-homogeneous elements $u,v\in V$ the half Jacobi identity holds:
 \begin{equation}\label{Hjacobi}
\begin{array}{c}
\displaystyle{z^{-1}_0\delta\left(\frac{z_1-z_2}{z_0}\right)
Y(u,z_1)Y(v,z_2)-\e _{u,v}\, z^{-1}_0\delta\left(\frac{z_2-z_1}{-z_0}\right)
Y(v,z_2)Y(u,z_1)}\\
\displaystyle{\simeq z_2^{-1}\delta\left(\frac{z_1-z_0}{z_2}\right)
Y(Y(u,z_0)v,z_2)}.
\end{array}
\end{equation}
Finally, for any $\Z _2$-homogeneous $u,v\in U$ and $n\geq 0$ we assume
\begin{equation}\label{Hhomog}
    |u_nv|=|u|+|v|
\end{equation}
(i.e. $u_n v$ is homogeneous and $|u_nv|=|u|+|v|$).

Sometimes we will emphasize that $Y(u,z)$ is pertinent to a vertex Lie
superalgebra $U$ by writing $Y_U(u,z)$.

Loosely speaking, a vertex Lie superalgebra, VLSA for short, carries the
``whole half''
of the structure of vertex superalgebra related to ``positive''
multiplications, except that properties of $\1$ in the definition of VSA are
replaced
by the half skew symmetry in the definition of VLSA. The half skew symmetry
can be
written by components as
\begin{equation}\label{HSScomp}
u_nv=-\e_{u,v}\,\sum_{k\geq 0}(-1)^{n+k}(D^k/k!)\,v_{n+k}u
\ \ \ \mbox{for all}\ \ n\geq 0.
\end{equation}
Also note that in the case when $U=U^0$, i.e. when all vectors are even, all
$\Z_2$-grading conditions become trivial and we speak of a vertex Lie
algebra,
VLA for short.

We define a homomorphism $\varphi : U\rightarrow W$ of two VLSA as a
$\Z_2$-grading
preserving linear map $\varphi$ such that
$$
\varphi (u_n v)=(\varphi (u))_n(\varphi (v)), \ \ \ \varphi D=D\varphi.
$$

Left (resp. right, two-sided) ideals in $U$ are defined as left (resp.
right, two-sided)
ideals for all multiplications. Note that an one-sided $\Z_2$-graded ideal
in VLSA
which is invariant
for $D$ must be two-sided ideal because of the half skew symmetry.
For a subset $S\subset U$ we denote by $\langle S \rangle$ a vertex Lie
superalgebra
generated by $S$, i.e. the smallest vertex Lie superalgebra containing the
set $S$.

A partial justification for our terminology might be the following
lemma proved in [B] by R. E. Borcherds in the case when $U$ is a
vertex algebra:
\begin{lem}\label{U/DU} Let $U$ be a VLSA. Then $U/DU$ is a Lie superalgebra
with
a commutator $[u+DU, v+DU]=u_0v+DU$.
\end{lem}
\pf Since $(Du)_0=0$, $DU$ is invariant for right multiplications ${}_0v$,
and then
by the skew symmetry (\ref{HSScomp}) for left multiplications $u_0$ as well.
Hence
on the quotient $U/DU$ we have a well defined bilinear operation $u_0v$ such
that $u_0v=-\e_{u,v}\,v_0u$. Now for $k=n=m=0$ in (\ref{jacobicomp}) we get
$u_0(v_0w)- \e_{u,v}v_0(u_0w)=(u_0v)_0w$, a Lie superalgebra Jacobi
identity. \qed

It is clear that any vertex superalgebra $V$ may be viewed as a VLSA.
Moreover, any
subspace $U\subset V$ invariant for $D$ and closed for multiplications
$u_nv$ for
$n\geq 0$ is a VLSA. Note that in this case
$$
Y_U(u,z)=Y^+_V(u,z)=\sum_{n\geq 0} u_nz^{-n-1}\simeq Y_V(u,z),
$$
so that $Y$ for $U$ should not be confused with $Y$ for $V$.

For a given vertex superalgebra $V$ we shall be especially interested
in vertex Lie superalgebras $\langle S \rangle\subset V$ generated by a
subspace
$S$ and the derivation $D$, i.e. when vertex Lie superalgebras
$\langle S \rangle\subset V$ are of the form
$$
\langle S \rangle=\sum_{k\geq 0}\,D^kS.
$$
In such a case any element $A$ in $\langle S \rangle$ can be written as a
combination
of vectors of the form $P(D)u$, $u\in S$, $P$ a polynomial, and (\ref{Hder})
implies
$$
Y_U(P(D)u,z)Q(D)v=P(d/dz)Q(D-d/dz)Y_U(u,z)v.
$$
Hence the structure of vertex Lie superalgebra $\langle S \rangle\subset V$
is
completely determined by $Y_S(u,z)v=Y_U(u,z)v$ for $u,v\in S$. Since
\begin{eqnarray*}
Y_V(u,z_1)v&=&\mbox{Res}_{z_2}\,z_2^{-1}Y_V(u,z_1)Y_V(v,z_2)\1\\
&\simeq& \mbox{Res}_{z_2}\,z_2^{-1}\left(Y_V(u,z_1)Y_V(v,z_2)\1
-\e_{u,v}\,Y_V(v,z_2)Y_V(u,z_1)\1\right)\\
&=&\mbox{Res}_{z_2}\,z_2^{-1}[Y_V(u,z_1),Y_V(v,z_2)]\1,
\end{eqnarray*}
we may use (for $m\geq 0$)
\begin{eqnarray*}
& & \mbox{Res}_{z_2}\,z_2^{-1}z_1^{-m-1}
\delta^{(m)}\left(\frac{z_2}{z_1}\right)e^{z_2D}w\simeq m!\,w\,z_1^{-m-1},\\
& & \mbox{Res}_{z_1}\,z_1^{-1}z_1^{-m-1}
\delta^{(m)}\left(\frac{z_2}{z_1}\right)e^{z_2D}w\simeq (-1)^m\, m!\,
\sum_{k=0}^m \frac{1}{k!}D^kw\,z_1^{-m-1+k},
\end{eqnarray*}
and easily calculate the principal part of the Laurent series
$Y_V(u,z)v$ from a given commutator. For example, in the case of
Neveu-Schwarz algebra (cf. [Li, (4.2.10)])
$$
[Y(\omega,z_1),Y(\tau,z_2)]=
z_1^{-1}\delta^{}\left(\frac{z_2}{z_1}\right)Y(D\tau,z_2)+
\frac{3}{2}z_1^{-2}\delta\,'\left(\frac{z_2}{z_1}\right)Y(\tau,z_2)
$$
gives
$$
Y(\omega,z)\tau \simeq
\sum_{n\geq 0}\omega_n\,\tau z^{-n-1}=\frac{D\tau}{z}+\frac{(3/2)\tau}{z^2}.
$$

The following examples of VLSA are obtained in this way from the well
known vertex (super)algebras associated to affine Lie
algebras, Virasoro algebra and Neveu-Schwarz algebra,
for details one may see [Li, Section 4]:

{\bf The case of affine Lie algebras.}
\ \ Here $S=S^0={\bf g} \oplus \C\1$ is a sum of
1-dimensional space and a Lie algebra {\bf g} with an invariant symmetric
bilinear
form $\langle\cdot,\cdot \rangle$,
$$
\langle S \rangle=\C\1\oplus\mbox{span}\,\{D^k x\mid x\in {\bf g}, k\geq 0\}
$$
and
$$
Y_S(x,z)y=\frac{[x,y]}{z}+\frac{\langle x,y\rangle \1}{z^2}, \ \ \ \
Y_S(x,z)\1=0, \ \ \ \
Y_S(\1,z)=0
$$
for $x,y\in {\bf g}$.

{\bf A case of affine Lie superalgebras.}
\ \ Here $S={\bf g} \oplus M\oplus \C\1$ is a sum of
1-dimensional space, a Lie algebra {\bf g} and a {\bf g}-module $M$,
where both {\bf g} and $M$ have {\bf g}-invariant symmetric bilinear forms
denoted as $\langle\cdot,\cdot \rangle$,
with $S^0={\bf g} \oplus \C\1$ and
$S^1=M $,
$$
\langle S \rangle=\C\1\oplus
\mbox{span}\,\{D^k x, D^k u\mid x\in {\bf g}, u\in M, k\geq 0\}
$$
and
\begin{eqnarray*}
&&Y_S(x,z)y=\frac{[x,y]}{z}+\frac{\langle x,y\rangle \1}{z^2}, \ \ \ \
Y_S(x,z)\1=0, \\
&&Y_S(u,z)v=\frac{\langle u,v\rangle \1}{z}, \ \ \ \
Y_S(u,z)\1=0, \\
&&Y_S(x,z)u=\frac{x\cdot u}{z}, \ \ \ \ Y_S(u,z)x=-\frac{x\cdot u}{z}, \\
&&Y_S(\1,z)=0
\end{eqnarray*}
for $x,y\in {\bf g}$ and $u,v\in M$.

{\bf The case of Virasoro algebra.}\ \ Here $S=S^0=\C\omega \oplus \C\1$ is
a
2-dimensional space and $\ell\in\C$,
$$
\langle S \rangle=\C\1\oplus
\mbox{span}\,\{D^k \omega\mid  k\geq 0\}
$$
and
$$
Y_S(\omega,z)\omega
=\frac{D\omega}{z}+\frac{2\omega}{z^2}+\frac{(\ell/2)\1}{z^4}, \ \ \ \
Y_S(\omega,z)\1=0, \ \ \ \
Y_S(\1,z)=0.
$$

{\bf The case of Neveu-Schwarz algebra.}\ \ Here
$S=S^0=\C\omega \oplus\C \tau \oplus \C\1$ is a
3-dimensional space and $\ell\in\C$, with $S^0=\C\omega \oplus \C\1$ and
$S^1=\C \tau $,
$$
\langle S \rangle=\C\1\oplus
\mbox{span}\,\{D^k \omega, D^k\tau\mid  k\geq 0\}
$$
and
\begin{eqnarray*}
&&Y_S(\omega,z)\omega
=\frac{D\omega}{z}+\frac{2\omega}{z^2}+\frac{(\ell/2)\1}{z^4}, \ \ \ \
Y_S(\omega,z)\1=0, \\
&&Y_S(\tau,z)\tau =\frac{2\omega}{z}+\frac{(2\ell/3)\1}{z^3}, \ \ \ \
Y_S(\tau,z)\1=0,\\
&&Y_S(\omega,z)\tau =\frac{D\tau}{z}+\frac{(3/2)\tau}{z^2}, \ \ \ \
Y_S(\tau,z)\omega =\frac{(1/2)D\tau}{z}+\frac{(3/2)\tau}{z^2},\\
&&Y_S(\1,z)=0.
\end{eqnarray*}

\begin{rem}\label{1=c}
Since for a vertex Lie superalgebra there is no notion of a vacuum vector,
it would be better to denote the vector $\1$
in the above examples by some other letter. So later on we shall set in the
affine case $S^0={\bf g} \oplus \C c$ with
$$
Y_S(x,z)y=\frac{[x,y]}{z}+\frac{\langle x,y\rangle c}{z^2};
$$
in the Virasoro case $S^0=\C\omega \oplus \C c$ with
$$
Y_S(\omega,z)\omega
=\frac{D\omega}{z}+\frac{2\omega}{z^2}+\frac{(1/2)c}{z^4}.
$$
\end{rem}

In all of these examples the set $S\subset V$ generates the corresponding
vertex superalgebra $V$.
In particular, the vertex Lie superalgebra $U=\langle S \rangle\subset V$
generates the
corresponding vertex superalgebra $V$, and obviously this is the smallest
such vertex
Lie superalgebra. On the other
extreme, for a $\frac12\Z_+$-graded vertex operator superalgebra $V$ we have
$U=\langle S\rangle=\oplus_{k\geq 0}D^kS=V$ if we take $S$ to be the set of
all quasi-primary fields, i.e. $S=\{u\in V\mid L(1)u=0\}$.

\section{Lie algebras associated to vertex Lie algebras}

Let $U$ be a vertex Lie superalgebra with the derivation $D$ and
consider the affinization $U \otimes \C[q,q^{-1}]$ as a $\Z_2$-graded tensor
product
of $U$ with the even space $\C[q,q^{-1}]$.
For vectors in $U \otimes \C[q,q^{-1}]$ we shall write
$$
  u_n = u \otimes q^n
$$
when $u \in U$ and $n \in \Z $. Consider the quotient $\Z_2$-graded vector
space
$$
\L (U) = \left(U \otimes \C[q,q^{-1}]\right)/
\mbox{span}\, \{(Du)_n + n u_{n-1} \mid u \in U , n \in \Z \}.
$$
The image of the vector $u_n $ in the quotient space we again denote by
$u_n$.
Note that we have three meanings for $u_n$, originally for $n\geq 0$
as a linear operator on $U$, now also for all $n\in \Z$ as an element of
$U \otimes \C[q,q^{-1}] $ or $\L (U)$, but it will be clear from
the context which meaning we have in mind.

\begin{thm}\label{L(U)lie-alg} $\L (U)$ is a Lie superalgebra with the
commutator defined by
\begin{equation}\label{L(U)comm}
[u_n, v_p] = \sum_{i \geq 0} {n \choose i}(u_iv)_{n+p-i}
\end{equation}
for $u,v\in U$ and $n,p\in\Z$. Moreover, we have the relation
\begin{equation}\label{L(U)der}
(Du)_n = - n \, u_{n-1}
\end{equation}
for all $u\in U$, $n\in\Z$, and the map $D : \L(U)\to \L(U)$ defined by
$D(u_n)=(Du)_n$ is an even derivation of the Lie superalgebra $\L(U)$.
\end{thm}
In the case when $U$ is a vertex algebra this theorem is a special case of
Borcherds' Lemma \ref{U/DU} applied to the affinization of $U$, in such a
case
a tensor product of vertex algebras. Direct proofs are also known, see
[Li2, Proposition 2.2.3 and Remark 2.2.4], here we follow [P] where the
techniques introduced in [FLM] are used. By following [FF] we may call
$\L(U)$ the local algebra of the vertex Lie algebra $U$.
This notion of local algebra $\L(U)$ is a special case of a more general
notion
of local vertex Lie algebra over the base space $U$ introduced in [DLM].

It will be convenient to consider formal Laurent series (vertex operators)
\begin{equation}\label{L(U)Y}
Y(u,z) = \sum_{n \in \Z} u_n z^{-n-1}
\end{equation}
with coefficients in $\L(U)$. Then we can write relations  (\ref{L(U)comm})
and (\ref{L(U)der}) as
\begin{equation}\label{L(U)commY}
[Y(u,z_1), Y(v,z_2)] = \mbox{Res}_{z_0} z^{-1}_2 \delta \left(\frac{z_1 -
z_0}{z_2}
\right) Y(Y(u,z_0)v,z_2),
\end{equation}
\begin{equation}\label{L(U)derY}
Y(Du,z) = \frac{d}{dz} Y(u,z) .
\end{equation}
Note that in (\ref{L(U)commY}) there are three $Y$'s with coefficients in
$\L(U)$ and
$Y(u,z_0)$ defined as $\sum_{n\geq 0}\,u_nz_0^{-n-1}$ has coefficients in
$\mbox{End}\,(U)$.
Also note
that in the expansion of $\delta$-function there are only nonnegative powers
of $z_0$,
so when we calculate the residue $\mbox{Res}_{z_0}$ we need to know only the
principal part of $Y(u,z_0)$. For this reason the commutator formula
(\ref{L(U)commY})
written in components looks the same as in the vertex algebra case. This is
the
key observation used in the proof of the above theorem.

\pf First note that for vectors  $u_n = u \otimes q^n $ the
formula (\ref{L(U)comm}) defines a bilinear operation $[\cdot,\cdot]$
which makes $U \otimes \C[q,q^{-1}]$ a $\Z_2$-graded algebra.
Since $D$ is an even derivation of $U$, it is obvious from (\ref{L(U)comm})
that $D=D\otimes 1$ is an even derivation of $[\cdot,\cdot]$.
The relations (\ref{Hder}) for the derivation $D$ of VLSA $U$ imply
$$
[(Du)_n + n u_{n-1} , v_p] = 0 , \ \ \ [u_n, (Dv)_p + p v_{p-1}] \equiv 0 .
$$
Here the first relation follows from $(Du)_i=-iu_{i-1}$, for the second we
also
need $u_i(Dv)=D(u_iv)-(Du)_iv$ and the result is in
$\mbox{span}\, \{(Dw)_n + n w_{n-1}\}$.
Hence the bilinear operation on the quotient space $\L(U)$ is well defined.
It is clear $\L(U)$ is $\Z_2$-graded algebra and that $D$ is well defined
even derivation such that (\ref{L(U)der}) holds. It remains
to check the axioms of Lie superalgebra. Since
\begin{eqnarray}
[Y(u,z_1), Y(v,z_2)]
&=&  \mbox{Res}_{z_0}\, z^{-1}_2 \delta
\left(\frac{z_1-z_0}{z_2}\right) Y(Y(u,z_0) v,z_2) \label{L(U)ss1}\\
&=& \e_{u,v}\, \mbox{Res}_{z_0}\, z^{-1}_2 \delta
\left(\frac{z_1-z_0}{z_2}\right) Y(e^{z_0D}Y(v,-z_0) u,z_2)
\label{L(U)ss2}\\
&=&  \e_{u,v}\,\mbox{Res}_{z_0}\, z^{-1}_1 \delta
\left(\frac{z_2+z_0}{z_1}\right) Y(Y(v,-z_0) u,z_2+ z_0) \label{L(U)ss3}\\
&=&  \e_{u,v}\,\mbox{Res}_{z_0}\, z^{-1}_1 \delta
\left(\frac{z_2+z_0}{z_1}\right) Y(Y(v,-z_0) u,z_1) \nonumber\\
&=&  - \e_{u,v}\,[Y (v,z_2), Y(u,z_1)],\nonumber
\end{eqnarray}
the skew symmetry for Lie superalgebra holds. Note that (\ref{L(U)ss2})
follows from (\ref{L(U)ss1}) by using the half skew symmetry (\ref{HSS}).
As it was already mentioned, for getting this equality it was sufficient to
know
the equality for principal parts of Laurent series in variable $z_0$.
Note that (\ref{L(U)ss3}) follows from (\ref{L(U)ss2}) because
(\ref{L(U)derY}) holds.

To prove the Jacobi identity for Lie superalgebras we use the definition of
commutator and get
\begin{eqnarray*}
& & [[Y(u,z_1), Y(v,z_2)],
Y(w,z_3)] \\
& = &
\mbox{Res}_{z_{23}} z^{-1}_2 \delta \left(\frac{z_3+z_{23}}{z_2}\right)
\mbox{Res}_{z_{12}} z^{-1}_1 \delta\left(\frac{z_2+z_{12}}{z_1}\right)
Y(Y(Y(u,z_{12}) v,z_{23})w,z_3)  .
\end{eqnarray*}
Note that in the expansions of $\delta$-functions there are only the
nonnegative powers of $z_{12}$ and $z_{23}$ and that, by taking the
residues, the above expression involves only the coefficients of
the principal part of $Y(Y(u,z_{12}) v,z_{23})w$. Hence we may apply
the half associator formula (\ref{Hassociator}) and get
\begin{eqnarray}
& & [[Y(u,z_1), Y(v,z_2)],Y(w,z_3)] \nonumber\\
& = &
\mbox{Res}_{z_{23}} z^{-1}_2 \delta \left(\frac{z_3+z_{23}}{z_2}\right)
\mbox{Res}_{z_{12}} z^{-1}_1 \delta\left(\frac{z_2+z_{12}}{z_1}\right)
\mbox{Res}_{z_{13}}z^{-1}_{12}\delta\left(\frac{z_{13}-z_{23}}{z_{12}}\right
)\nonumber \\
&&\times Y\left(Y(u,z_{13})Y(v,z_{23})w,z_3\right)\label{L(U)jacobi1} \\
&&-\e _{u,v}\,
\mbox{Res}_{z_{23}} z^{-1}_2 \delta \left(\frac{z_3+z_{23}}{z_2}\right)
\mbox{Res}_{z_{12}} z^{-1}_1 \delta\left(\frac{z_2+z_{12}}{z_1}\right)
\mbox{Res}_{z_{13}}
z^{-1}_{12}\delta\left(\frac{z_{23}-z_{13}}{-z_{12}}\right)\nonumber \\
&&\times Y\left(Y(v,z_{23})Y(u,z_{13})w,z_3\right) \label{L(U)jacobi2}  \\
& = &
\mbox{Res}_{z_{13}}\mbox{Res}_{z_{23}}
z^{-1}_1 \delta\left(\frac{z_3+z_{13}}{z_1}\right)
z^{-1}_2 \delta \left(\frac{z_3+z_{23}}{z_2}\right) \nonumber \\
&&\times Y\left(Y(u,z_{13})Y(v,z_{23})w,z_3\right)\label{L(U)jacobi3} \\
&&-\e _{u,v}\, \mbox{Res}_{z_{23}}\mbox{Res}_{z_{13}}
z^{-1}_1 \delta\left(\frac{z_3+z_{13}}{z_1}\right)
 z^{-1}_2 \delta \left(\frac{z_3+z_{23}}{z_2}\right) \nonumber \\
&&\times Y\left(Y(v,z_{23})Y(u,z_{13})w,z_3\right) \label{L(U)jacobi4}\\
&=&
[Y(u,z_1), [Y(v,z_2),Y(w,z_3)]] -\e _{u,v}\,[Y(v,z_2),[Y(u,z_1),Y(w,z_3)]],
\nonumber
\end{eqnarray}
which is the Jacobi identity for Lie superalgebra $\L(U)$. Note that we
obtained the expression (\ref{L(U)jacobi3}) from (\ref{L(U)jacobi1}) by
using the
identity for $\delta$-functions:
\begin{eqnarray}
&& \mbox{Res}_{z_{12}}
z^{-1}_2 \delta \left(\frac{z_3+z_{23}}{z_2}\right)
z^{-1}_1 \delta\left(\frac{z_2+z_{12}}{z_1}\right)
z^{-1}_{12}\delta\left(\frac{z_{13}-z_{23}}{z_{12}}\right)\nonumber \\
& = &
\mbox{Res}_{z_{12}}
z^{-1}_2 \delta \left(\frac{z_3+z_{23}}{z_2}\right)
z^{-1}_{12}\delta\left(\frac{z_{13}-z_{23}}{z_{12}}\right)
\sum_n\sum_{r\geq 0}{n \choose r} z_2^{n-r} z_{12}^r z_1^{-n-1}\nonumber \\
& = &
\mbox{Res}_{z_{12}}
z^{-1}_2 \delta \left(\frac{z_3+z_{23}}{z_2}\right)
z^{-1}_{12}\delta\left(\frac{z_{13}-z_{23}}{z_{12}}\right)\nonumber \\
&&\sum_n\sum_{r\geq 0}{n \choose r}
(z_3+z_{23})^{n-r} (z_{13}-z_{23})^r z_1^{-n-1}\nonumber \\
& = &
z^{-1}_2 \delta \left(\frac{z_3+z_{23}}{z_2}\right)
\sum_n\sum_{r\geq 0}\sum_{s\geq 0}{n \choose r} {{n-r} \choose s}
z_3^{n-r-s}z_{23}^s (z_{13}-z_{23})^r z_1^{-n-1}\nonumber \\
& = &
z^{-1}_2 \delta \left(\frac{z_3+z_{23}}{z_2}\right)
\sum_n\sum_{k\geq 0}\sum_{s+r=k}{n \choose k} {k \choose r}
z_3^{n-k}z_{23}^s (z_{13}-z_{23})^r z_1^{-n-1}\nonumber \\
& = &
z^{-1}_2 \delta \left(\frac{z_3+z_{23}}{z_2}\right)
\sum_n\sum_{k\geq 0}{n \choose k}
z_3^{n-k} z_{13}^k z_1^{-n-1}\nonumber \\
& = &
z^{-1}_2 \delta \left(\frac{z_3+z_{23}}{z_2}\right)
z^{-1}_1 \delta\left(\frac{z_3+z_{13}}{z_1}\right). \nonumber
\end{eqnarray}
This identity for $\delta$-functions should be understood in the context
of the above proof, i.e. in the presence of the term
$Y\left(Y(u,z_{13})Y(v,z_{23})w,z_3\right)$ and the residues
$\mbox{Res}_{z_{13}}\mbox{Res}_{z_{23}}$.
We obtain the expression (\ref{L(U)jacobi4}) from (\ref{L(U)jacobi2})
in a similar way. \qed

\begin{rem}\label{defDLM}
Motivated by a definition of local vertex Lie algebra over the base space
$U$
given in [DLM] and by methods used in [KW], let us note the following:

Given a vector space $U$ and a linear map $Y$ satisfying (\ref{Hvlw0}),
(\ref{Hder}) and
(\ref{Hhomog}), later on we shall say that $U$ is a VLA-J-SS, we have on the
quotient
space $\L(U)$ a bilinear operation $[\cdot,\cdot]$ defined by
(\ref{L(U)comm}).
Then $\L(U)$ is a Lie (super)algebra if and only if $U$ is a vertex Lie
(super)algebra.

To see the ``only if'' part, note that, by Theorem \ref{UinL(U)} below,
$u=0$ if and only
if $Y(u,z)=0$. The skew symmetry for $[\cdot,\cdot]$ implies that
(\ref{L(U)ss1})=(\ref{L(U)ss2}) and, by the above remark and Lemma 2.1.4 in
[Li], this
implies the half skew symmetry (\ref{HSS}). In a similar way we see that the
Jacobi
identity for $[\cdot,\cdot]$ implies the half associator formula
(\ref{Hassociator}), which by Lemma \ref{assoc} below implies the half
Jacobi identity
(\ref{Hjacobi}).
\end{rem}

\begin{prop}\label{L(U)toL(U)}
Let $\varphi : U\to W$ be a homomorphism of vertex Lie superalgebras. Then
\begin{equation}\label{L(phi)}
\L(\varphi) : \L(U) \to \L(W), \ \ \ u_n \mapsto \left(\varphi (u)\right)_n
\end{equation}
is a homomorphism of Lie superalgebras. Moreover,
$\L(\varphi)D=D\L(\varphi)$.
\end{prop}

\pf Since $\varphi D=D\varphi$, the map $u\otimes q^n \to \varphi (u)\otimes
q^n$
factors through the quotient and the map $\L(\varphi)$ is well defined
$\Z_2$-grading preserving map such that $\L(\varphi)D=D\L(\varphi)$.
From the definition of commutator (\ref{L(U)comm}) we have
$$
 \sum_{i \geq 0} {n \choose i}(u_iv)_{n+p-i}\  \mapsto \
 \sum_{i \geq 0} {n \choose i}\left(\varphi (u_iv)\right)_{n+p-i}
= \sum_{i \geq 0} {n \choose i}\left(\varphi (u)_i
\,\varphi(v)\right)_{n+p-i}\ ,
$$
i.e. $\L(\varphi)$ is a homomorphism of Lie superalgebras. \qed
$$
\mbox{Set} \ \ \ \ \L_-(U)=\mbox{span}\,\{u_n\mid u\in U,\, n<0\}, \ \ \
\L_+(U)=\mbox{span}\,\{u_n\mid u\in U,\, n\geq 0\}.
$$
\begin{prop}\label{L(U)pm}
$\L_-(U)$ and $\L_+(U)$ are Lie superalgebras invariant for $D$ and
\begin{eqnarray}
&&\L(U)=\L_-(U)\oplus \L_+(U),\nonumber\\
&&\L_-(U)\cong \left(U\otimes q^{-1}\C[q^{-1}]\right)/
\mbox{span}\, \{(Du)_n + n u_{n-1} \mid u \in U , n < 0\},\nonumber\\
&&\L_+(U)\cong \left(U\otimes \C[q]\right)/
\mbox{span}\, \{(Du)_n + n u_{n-1} \mid u \in U , n \geq 0\}.\nonumber
\end{eqnarray}
\end{prop}

\pf It is clear from the definition (\ref{L(U)comm}) of commutator in
$\L(U)$ that both $\L_-(U)$ and $\L_+(U)$ are subalgebras. The invariance
for $D$ follows from (\ref{L(U)der}). The other statements follow from the
fact that $\L(U)$ is the quotient of
$$
\left(U\otimes q^{-1}\C[q^{-1}]\right)\oplus \left(U\otimes \C[q]\right)
$$
by the sum of $\Z_2$-graded $D$ invariant subspaces
$$
\mbox{span}\, \{(Du)_n + n u_{n-1} \mid u \in U , n < 0\}\oplus
\mbox{span}\, \{(Du)_n + n u_{n-1} \mid u \in U , n \geq 0 \}.\ \ \Box
$$

The following proposition is an obvious consequence of the definition
(\ref{L(phi)}) of $\L(\varphi)$:

\begin{prop}\label{L(U)pmtoL(U)pm}
Let $\varphi : U\to W$ be a homomorphism of vertex Lie superalgebras and let
$\L_\pm(\varphi)$ be the restriction of $\L(\varphi)$ on $\L_\pm(U)$. Then
$\L_\pm(\varphi)(\L_\pm(U))\subset \L_\pm(W)$ and
$$
\L_-(\varphi) : \L_-(U) \to \L_-(W), \ \ \ \L_+(\varphi) : \L_+(U) \to
\L_+(W)
$$
are homomorphisms of Lie superalgebras. Moreover,
$\L_\pm(\varphi)D=D\L_\pm(\varphi)$.
\end{prop}

\begin{thm}\label{UinL(U)}
The map
$$
\iota_U : U\to \L_-(U), \ \ \ u\mapsto u_{-1}
$$
is an isomorphism of $\Z_2$-graded vector spaces $U$ and $\L_-(U)$ and
$\iota_U D=D\iota_U$.

Moreover, if  $\varphi : U\to W$ is a homomorphism of vertex Lie
superalgebras, then
$\iota_W \varphi=\L_-(\varphi)\iota_U$.
\end{thm}

\pf It follows from (\ref{L(U)der}) that
$u_{-k-1}=\frac{1}{k!}\,(D^ku)_{-1}$
for $k\geq 0$, and this implies that $\iota_U$ is surjective. Let
$u_{-1}=0$, that is
$$
u\otimes q^{-1}=\sum_{i=-n}^m \left(Dv^{(i)}\otimes q^{i}+iv^{(i)}\otimes
q^{i-1}\right)
$$
for some $n,m\geq 1$ and some $v^{(i)}\in U$. Obviously
$\sum_{i\geq 0}( Dv^{(i)}\otimes q^{i}+iv^{(i)}\otimes q^{i-1})=0$.
Since $v^{(i)}=0$ implies $Dv^{(i)}=0$, by induction we get
$v^{(i)}=0$ for all $i< 0$, and hence $u=0$.
So $\iota_U$ is an isomorphism. The remaining statements are clear. \qed

We may identify $U$ with $\L_-(U)$ via the map $\iota_U$ and consider
$$
U\subset \L(U)=U\oplus \L_+(U).
$$
If we transport the Lie superalgebra structure of $\L_-(U)$ on $U$, then the
commutator (\ref{L(U)comm}) reads for $u,v\in U$ as
\begin{equation}\label{Ucomm}
[u,v]=\sum_{n\geq 0}\ ((-1)^n/(n+1)!)\, D^{n+1}\,(u_nv).
\end{equation}
Note that by Theorem \ref{UinL(U)} any homomorphism
$\varphi : U\to W$ of vertex Lie superalgebras
is a homomorphism of Lie superalgebras $U$ and $W$ with commutators defined
by
(\ref{Ucomm}). Moreover, this homomorphism $\varphi$ extends to
a homomorphism of Lie superalgebras $\L(U)$ and $\L(W)$.


\section{Enveloping vertex algebras of vertex Lie algebras}

Let $U$ be a vertex Lie superalgebra with the derivation $D$ and let
$\L(U)=\L_-(U)\oplus \L_+(U)$ be the corresponding  Lie superalgebra
with the derivation $D$. The induction by a
trivial $\L_+(U)$-module $\C$ gives a generalized Verma $\L(U)$-module
\[
   {\cal V}(U)={\cal U}\left({\cal L}(U)\right)
\otimes_{\,{\cal U}\left({\cal L}_+(U)\right)} \C ,
\]
where $\cal U$ stands for the universal enveloping algebra of a given Lie
superalgebra.
The $\L(U)$-module $\V(U)$ is isomorphic to a quotient
$$
\V(U)=\U(\L(U))/\U(\L(U))\L_+(U)
$$
of $\U(\L(U))$ by a left ideal generated by $\L_+(U)$.  Clearly $\V(U)$ is a
$\Z_2$-graded space and the action of Lie superalgebra $\L(U)$ is given by
the left multiplication. Note that the derivation $D$ of Lie superalgebra
$\L(U)$ extends to
a derivation $D$ of the associative superalgebra $\U(\L(U))$, and since $D$
preserves $\L_+(U)$, it defines an even operator $D$ on the quotient
$\V(U)$.
We denote by $\1\in\V(U)$ the image of $1\in \U(\L(U))$. It is an even
vector and we have
$$
D\1=0.
$$
If we think of $\V(U)$ only as a $\Z_2$-graded space, then
$$
   {\cal V}(U)\cong {\cal U}\left({\cal L}(U)_-\right)
$$
and under this identification $\1$ is the identity and $D$ is an even
derivation of the associative superalgebra $\U(\L_-(U))$ extending the
derivation $D$
of $\L_-(U)$.
By Theorem \ref{UinL(U)} we may identify $U$ and $\L_-(U)$ via the map
$\iota_U$,
so clearly the map
$$
\kappa_U : U\to \V(U),\ \ \ \ u\mapsto u_{-1}\1
$$
is an injection. Sometimes it will be convenient to identify
$U$ and the $\Z_2$-graded subspace $\kappa_U (U)\subset \V(U)$, i.e. to
consider
\begin{equation}\label{UinV(U)2}
U\subset \V(U).
\end{equation}
\begin{lem}\label{V(U)der1}
For the action of $u_n\in\L(U)$ on the module $\V(U)$ we have
\begin{equation}\label{V(U)der2}
[D,u_n]=(Du)_n=-n u_{n-1}.
\end{equation}
\end{lem}
\pf By definition the operator $D$ acts on a vector
\begin{equation}\label{vectinV(U)}
w=u^{(1)}\dots u^{(k)}\1\in\V(U), \ \ \ \ u^{(1)},\dots ,u^{(k)}\in\L(U),
\end{equation}
as a derivation, and by definition $u_n$ acts by the left multiplication, so
$Du_nw-u_nDw=D(u_n)w$. Hence $[D,u_n]=D(u_n)=(Du)_n$ and the relation
(\ref{L(U)der})
implies the lemma. \qed

Since $D\1=0$, we have $D(u_{-1}\1)=[D,u_{-1}]\1=(Du)_{-1}\1$. Hence
$D\kappa_U=\kappa_U D$ and on $U\subset \V(U)$ both derivations coincide.

Since $\V(U)$ is a $\L(U)$-module, the formal Laurent
series $Y(u,z)$ defined by (\ref{L(U)Y}) operate on $\V(U)$ and the
corresponding formal Laurent series we shall denote by $Y_{\V(U)}(u,z)$.
Hence the coefficients $u_n$ in
$$
Y_{\V(U)}(u,z) = Y_{\V(U)}(u_{-1}\1,z) = \sum_{n \in \Z} u_n z^{-n-1} \, ,\
\ \ \ u\in U
$$
are operators on $\V(U)$. Then we can restate (\ref{V(U)der2}) as
\begin{equation}\label{V(U)derY}
[D,Y_{\V(U)}(u,z)] =Y_{\V(U)}(Du,z) =\frac{d}{dz}\,Y_{\V(U)}(u,z).
\end{equation}
\begin{lem}\label{vertop} The set of formal Laurent series
$\{Y_{\V(U)}(u,z)\mid u\in U^0\cup U^1\}$ is a set of mutually local vertex
operators on $\V(U)$.
\end{lem}
\pf We have already proved (\ref{V(U)derY}), that is (\ref{derv(z)}) in
the definition of vertex operators. By construction we have
\begin{equation}\label{vn10V(U)}
u_n\1=0 \ \ \ \ \mbox{for all}\ \ u\in U,\  n\geq 0.
\end{equation}
Hence for a vector $w\in \V(U)$ of the form  (\ref{vectinV(U)}) we have
$$
u_nw=u_nu^{(1)}\dots u^{(k)}\1=[u_n,u^{(1)}\dots u^{(k)}]\1
$$
for $n\geq 0$, and then the commutator formula (\ref{L(U)comm}) implies
$u_nw=0$ for
sufficiently large $n$, i.e. (\ref{vlw0M}) in
the definition of vertex operators. The commutator formula (\ref{L(U)commY})
for $u,v\in U^0\cup U^1$ written in the form (\ref{commutator}) clearly
implies locality and the lemma follows. \qed

Since (\ref{vn10V(U)}) means that $Y_{\V(U)}(u,z)\1$ is a power series in
$z$ and since
by construction $\lim_{z\to 0}Y_{\V(U)}(u,z)\1 =u_{-1}\1=u$, all the
assumptions of
Theorem \ref{MPThm} hold and we have the following:
\begin{thm}\label{vaV(U)}
The set $\{Y_{\V(U)}(u,z)\mid u\in U^0\cup U^1\}$ of mutually local vertex
operators on $\V(U)$ generates the vertex superalgebra structure on $\V(U)$
with the vacuum vector $\1$ and the derivation $D$.
\end{thm}
\begin{prop}\label{UinV(U)}
The map $\kappa_U : U\to \V(U)$ is an injective homomorphism of
vertex Lie superalgebras.
\end{prop}
\pf
We have already seen that $\kappa_U$ is an injective $\Z_2$-grading
preserving map
such that $\kappa_U D=D\kappa_U$. So let $n\geq 0$. Because of
(\ref{vn10V(U)})
and the commutator formula (\ref{L(U)comm}) we have
$$
u_n(v_{-1}\1)=[u_n,v_{-1}]\1=\sum_{i \geq 0} {n \choose i}(u_iv)_{n-1-i}\1
={n \choose n}(u_nv)_{-1}\1.
$$
Hence
$$
\kappa_U (u_nv)=(u_nv)_{-1}\1=u_n(v_{-1}\1)=(u_{-1}\1)_n(v_{-1}\1)=
(\kappa_U (u))_n(\kappa_U (v)).\ \ \
\Box$$

The above proposition means that, with the identification (\ref{UinV(U)2}),
the
structure of vertex superalgebra on $\V(U)$  extends the structure of
vertex Lie superalgebra on $U\subset \V(U)$. As suggested by the next
theorem,
we shall say that $\V(U)$ is the universal enveloping vertex superalgebra of
the vertex Lie superalgebra $U$.

\begin{thm}\label{ueva}
Let $V$ be a vertex superalgebra and $\varphi : U\to V$ a homomorphism
of vertex Lie superalgebras. Then $\varphi$ extends uniquely to a
vertex superalgebra homomorphism $\tilde\varphi : \V(U)\to V$.
\end{thm}

\pf
Set $M=\V(U)\oplus V$ and for $u\in U$ define
$$
Y_{M}(u,z)=Y_{\V(U)}(u,z)+Y_{V}(\varphi (u),z).
$$
Then $\{Y_{M}(u,z)\mid u\in U^0\cup U^1\}$ is a set of mutually local vertex
operators on $M$ which, by Li's Theorem \ref{LiThm}, generates a vertex
algebra $W$.
Consider the restriction maps
$$
\V(U)\stackrel{p_1}{\longleftarrow} W \stackrel{p_2}{\longrightarrow} V
$$
defined by
\begin{eqnarray*}
&&p_1 : a(z)\mapsto a(z)|\V(U)\mapsto \lim_{z\to 0} \,(a(z)|\V(U))\1, \\
&&p_2 : a(z)\mapsto a(z)| V \mapsto \lim_{z\to 0} \,(a(z)|V)\1.
\end{eqnarray*}

\noindent{\bf Step 1.}\ \ {\it The maps $p_1$ and $p_2$ are
homomorphisms of vertex superalgebras.}
Recall that for two formal Laurent series $u(z), v(z)\in F(M)$ the product
$u(z)_nv(z)$ is defined by (\ref{product}), so it is clear that for
any invariant subspace $N\subset M$ we have
$$
(u(z)_nv(z))|N=(u(z)|N)_n(v(z)|N).
$$
For this reason both $ a(z)\mapsto a(z)|\V(U)$ and $a(z)\mapsto a(z)| V$
are homomorphisms of vertex superalgebras, the first restriction
from $W$ to a vertex superalgebra $W_1$ of vertex operators on $\V(U)$
generated by
$(Y_{M}(u,z))|\V(U)=Y_{\V(U)}(u,z)$, $u\in U$, the second restriction
from $W$ to a vertex superalgebra $W_2$ of vertex operators on $V$ generated
by
$(Y_{M}(u,z))|V=Y_{V}(\varphi(u),z)$, $u\in U$. Since vertex superalgebras
are
closed for multiplications, we obviously have
$$
W_1\subset \{Y_{\V(U)}(v,z)\mid v\in \V(U)\},\ \ \ \ \
W_2\subset \{Y_{V}(v,z)\mid v\in V\}.
$$
Since in general for any vertex superalgebra the map
$$
Y(v,z)\mapsto v_{-1}\1=\lim_{z\to 0} \,Y(v,z)\1
$$
is an isomorphism of vertex superalgebra of fields with the algebra itself,
both
$p_1$ and $p_2$ are homomorphisms of vertex superalgebras.

Let us note at this point that $p_1$ is a surjection since
$\V(U)\cong \{Y_{\V(U)}(v,z)\mid v\in \V(U)\}$
is generated by $\{Y_{\V(U)}(u,z)\mid u\in U\}$.

\noindent{\bf Step 2.}\ \ {\it $W$ is a $\L(U)$-module and $W=\U(\L_-(U))
I(z)$.} \
 A linear map
$$
U\otimes \C[q,q^{-1}]\to \mbox{End}\, W,\ \ \
u\otimes q^n\mapsto u(z)_n=Y_M(u,z)_n
$$
is obviously well defined. Since
\begin{eqnarray*}
 Du\otimes q^n+ n\,u\otimes q^{n-1}&\mapsto &(Du)(z)_n + n \,u(z)_{n-1}\\
&=& Y_{\V(U)}(Du,z)_n+Y_{V}(\varphi (Du),z)_n\\
& & +n\,Y_{\V(U)}(u,z)_{n-1}+n\,Y_{V}(\varphi (u),z)_{n-1}\\
&=& (\frac{d}{dz}\,Y_{\V(U)}(u,z))_n
  +(\frac{d}{dz}\,Y_{V}(\varphi (u),z)_n\\
& & +n\,Y_{\V(U)}(u,z)_{n-1}+n\,Y_{V}(\varphi (u),z)_{n-1}\\
&=&(\frac{d}{dz}\,u(z))_n + n \,u(z)_{n-1}=0,
\end{eqnarray*}
we have a well defined map on the quotient
$$
\L(U)\to \mbox{End}\, W,\ \ \ \ \
u_n\mapsto u(z)_n=Y_M(u,z)_n \ .
$$
Since $\varphi : U\to V$ is VLA homomorphism, we have for $n\geq 0$
\begin{eqnarray*}
(u_nv)(z)
&=& Y_{\V(U)}(u_nv,z)+Y_{V}(\varphi (u_nv),z)\\
&=& Y_{\V(U)}(u_nv,z)+Y_{V}(\varphi (u)_n\,\varphi (v),z)\\
&=& Y_{\V(U)}(u,z)_n\,Y_{\V(U)}(v,z)+Y_{V}(\varphi (u),z)_n\,Y_{V}(\varphi
(v),z)\\
&= &u(z)_nv(z),
\end{eqnarray*}
which together with (\ref{L(U)comm}) and (\ref{commutator})
implies
\begin{eqnarray}
[u_p, v_q] = \sum_{n \geq 0} {p \choose n}(u_nv)_{p+q-n} &\mapsto&
\sum_{n \geq 0} {p \choose n}(u_nv)(z)_{p+q-n}\nonumber\\
&=&\sum_{n \geq 0} {p \choose n}(u(z)_nv(z))_{p+q-n}\nonumber\\
&=&[u(z)_p, v(z)_q].\nonumber
\end{eqnarray}
Hence the map $u_n\mapsto u(z)_n$ is a representation of $\L(U)$ on $W$.

By definition $W$ is generated by $\1=I(z)$ and the set of homogeneous
fields $u(z)=Y_M(u,z)$, $u\in U^0\cup U^1$. This means that we consider all
possible products like
$$
((u^{(1)}(z)_{n_1}u^{(2)}(z))_{n_2}(u^{(3)}(z)_{n_3}
(u^{(4)}(z)_{n_4}u^{(5)}(z))).
$$
It is easy to see by using the associator formula
(cf. (\ref{jacobicomp})) that $W$ is spanned by elements of the
form
\begin{eqnarray}
&&u^{(1)}(z)_{n_1}(u^{(2)}(z)_{n_2}(\dots
(u^{(k)}(z)_{n_k}u^{(k+1)}(z))\dots ))\nonumber\\
&=&u^{(1)}(z)_{n_1}u^{(2)}(z)_{n_2}\dots
u^{(k)}(z)_{n_k}u^{(k+1)}(z)_{-1}\,I(z),
\label{monomial}
\end{eqnarray}
$u^{(1)}, u^{(2)},\dots ,u^{(k+1)}\in U$, $n_1, n_2,\dots ,n_{k}\in\Z$.
Hence $W=\U(\L(U)) I(z)$. Since we have $u(z)_nI(z)=0$ for $n\geq 0$, by
using the PBW theorem we finally get $W=\U(\L_-(U)) I(z)$.

\noindent{\bf Step 3.}\ \ {\it The map $p_1 : W\to \V(U)$ is an
isomorphism of vertex superalgebras.} \
We already know that $p_1$ is a surjective homomorphism. If we are to
identify $v\leftrightarrow Y_{\V(U)}(v,z)$ for $v\in \V(U)$, the map $p_1$
is just a restriction $ a(z)\mapsto a(z)|\V(U)$, so on monomials
of the form (\ref{monomial}) we have with our identification
$$
p_1 : u^{(1)}(z)_{n_1}u^{(2)}(z)_{n_2}\dots
u^{(k)}(z)_{n_k}u^{(k+1)}(z)_{n_{k+1}}\,I(z)
\mapsto
u^{(1)}_{n_1}u^{(2)}_{n_2}\dots u^{(k)}_{n_k}u^{(k+1)}_{n_{k+1}}\1.
$$
Hence for $u\in \U(\L_-(U))$ we have $p_1(u\,I(z))=u\1\cong u\in
\U(\L_-(U))$.
This implies that $p_1$ is an injection as well and that we have a
homomorphism of vertex superalgebras
$$
 \V(U)\stackrel{p_1^{-1}}{\longrightarrow} W \stackrel{p_2}{\longrightarrow}
V.
$$
Set $\tilde\varphi=p_2\circ p_1^{-1}$. Then by construction we have for
$u\in U$
$$
u\mapsto Y_M(u,z) \mapsto Y_V(\varphi(u),z) \leftrightarrow \varphi(u).
$$
Hence $\tilde\varphi|U=\varphi$. Since $U$ generates the vertex superalgebra
$\V(U)$,
the homomorphism $\tilde\varphi$ is uniquely determined by $\varphi$. \qed

\begin{coro}\label{V(phi)}
Let  $\varphi : U_1\to U_2$ be a homomorphism
of vertex Lie superalgebras $U_1$ and $U_2$. Then $\varphi$ extends uniquely
to a
vertex superalgebra homomorphism $\V(\varphi) : \V(U_1)\to \V(U_2)$.
\end{coro}

We shall say that a $\L(U)$-module $M$ is restricted if for any $u\in U$,
$w\in M$
$$
u_nw=0\ \ \      \mbox{for}\ \ n \ \ \mbox{sufficiently\ large.}
$$
Let us denote by $\L(U)\times \C D$ a Lie superalgebra
with $[D,u_n]=D(u_n)$ and let us say that a $(\L(U)\times \C D)$-module $M$
is restricted if $M$ is restricted as an $\L(U)$-module. Note that
$\V(U)$ is a restricted $(\L(U)\times \C D)$-module.

If $M$ is a $\L(U)$-module, then for any $u\in U$ we can form a formal
Laurent series
$Y_M(u,z)=\sum_{n\in \Z}\,u_nz^{-n-1}$ with elements $u_n\in \L(U)$ acting
as operators on the module $M$. If $M$ is a restricted $(\L(U)\times \C
D)$-module,
then $\{Y_M(u,z)\mid u\in U^0\cup U^1\}$ is a set of mutually local vertex
operators on $M$.

\begin{lem}\label{restrmod}
Let $M$ be a restricted $(\L(U)\times \C D)$-module and let $W$ be a vertex
superalgebra generated by $\{Y_M(u,z)\mid u\in U^0\cup U^1\}$.
Set $u(z)=Y_M(u,z)$ and denote by
$u(z)_n$ operators on $W$ defined by multiplications (\ref{product}).
Then the linear map
$$
\L(U)\to \mbox{End}\,W,\ \ \ \ u_n\mapsto u(z)_n
$$
is well defined and
\begin{eqnarray*}
&&(u_nv)(z)=u(z)_nv(z) \ \ \ \ \ \mbox{for all}\ \ n\geq 0,\\
&&[u_n,v_m]\mapsto [u(z)_n,v(z)_m] \ \ \ \ \ \mbox{for all}\ \ n,m \in \Z\,.
\end{eqnarray*}
In another words, the map
\begin{equation}\label{restrmod2}
U\to W,\ \ \ \ u\mapsto u(z)
\end{equation}
is a homomorphism of vertex Lie superalgebras and $W$ is a restricted
$(\L(U)\times \C D)$-module.
\end{lem}

\pf
Set $N=\V(U)\oplus M$ and for $u\in U$
$$
Y_{N}(u,z)=Y_{\V(U)}(u,z)+Y_{M}(u,z).
$$
This vertex operators on $N$ generate a vertex superalgebra $V$, denote by
$p_1$ a restriction map from $V$ to $\V(U)$ and by
$p_2$ a restriction map from $V$ to $W$. Then both $p_1$ and $p_2$
are homomorphisms of vertex superalgebras and both $\V(U)$ and $W$ are
$V$-modules.
By [Li, Lemma 2.3.5] the ``structure constants'' in the commutator formula
for $V$ are completely determined by the ``structure constants'' appearing
in the commutator formula for a faithful $V$-module.
Since for $u(z)\in V$ we have $Y_{\V(U)}(u(z),z_1)=Y_{\V(U)}(u,z_1)$, i.e.
$u(z)_n$ acts on $\V(U)$ as $u_n$, the commutator
formula for $V$-module $\V(U)$ reads
$$
[u(z)_p, v(z)_q]=[u_p, v_q] = \sum_{n \geq 0} {p \choose n}(u_nv)(z)_{p+q-n}
$$
and implies the commutator formula for $V$-module $M$
$$
[u(z)_p, v(z)_q]=\sum_{n \geq 0} {p \choose n}(u_nv)(z)_{p+q-n}.
$$
Since $[u(z)_p, v(z)_q]$ equals
$\sum_{n \geq 0} {p \choose n}(u(z)_nv(z))_{p+q-n}$,
the lemma follows. \qed

Since by Theorem \ref{ueva} a VLSA homomorphism (\ref{restrmod2}) extends to
a
VSA homomorphism $\V(U)\to W$, and since by Li's Theorem \ref{LiThm} $M$ is
a
$W$-module, we have the following:

\begin{thm}\label{V(U)modules}
Any restricted $(\L(U)\times \C D)$-module is a $\V(U)$-module.
\end{thm}

\begin{rem}\label{dokazi}
The proofs of Theorem \ref{ueva}, Lemma \ref{restrmod} and Theorem
\ref{V(U)modules}
are modeled after the arguments in [Li, Section 4]. Also note that Theorem
\ref{vaV(U)}
can be proved without using Theorem \ref{MPThm}, but rather by using the
proof of
Theorem \ref{ueva} with Step 2 changed with the help of Lemma
\ref{restrmod}. This would
be parallel to the argument in [Li] when proving that
for an affine Lie algebra $\tilde{\bf g}$ the generalized Verma
module $M_{\bf g}(\ell,\C)$ is a vertex algebra and that any
restricted $\tilde{\bf g}$-module of level $\ell$ is a $M_{\bf
g}(\ell,\C)$-module.

The results in this section are similar to some results in [DLM].
\end{rem}


\section{Commutativity and the skew symmetry for VLA}

For the purposes of the last two sections let us make a few technical
definitions:
Let $U$ be a $\Z _2$-graded vector space equipped with an even linear
operator $D$ on $U$ called the derivation and a linear map
$U\rightarrow z^{-1}(\mbox{End}\,U)[[z^{-1}]]$,
$u\mapsto Y(u,z)=\sum_{n\geq 0}u_nz^{-n-1}$,
satisfying the following conditions for homogeneous $u,v\in U$:
\begin{eqnarray*}
&&u_nv=0\ \ \
\mbox{for}\ \ n \ \ \mbox{sufficiently\ large},\\
&&(Du)_nv=-nu_{n-1}v ,\\
&&D(u_nv)=(Du)_{n}v+u_n(Dv), \\
&& |u_nv|=|u|+|v|.
\end{eqnarray*}
Then we shall say that $U$ is a VLA-J-SS (something like a vertex Lie
superalgebra without the half Jacobi identity and without the half
skew symmetry).

If for an VLA-J-SS $U$ the half Jacobi identity (\ref{Hjacobi}) holds, then
we shall say
that $U$ is a VLA-SS.

If for an VLA-SS $U$ the half skew symmetry (\ref{HSS}) holds, then clearly
$U$ is a
vertex Lie superalgebra.

For $U$ and $W$ having any of the above structures
we define a homomorphism $\varphi$ to be a $\Z_2$-grading
preserving linear map $\varphi : U\rightarrow W$ such that
$$
\varphi (u_n v)=(\varphi (u))_n(\varphi (v)), \ \ \ \varphi D=D\varphi.
$$

Left (resp. right, two-sided) ideals in $U$ are defined as left (resp.
right, two-sided)
ideals for all multiplications.

\begin{lem}\label{comm} In the definition of VLA-SS the half Jacobi identity
can be
equivalently substituted by the half commutator formula.
\end{lem}

\pf Let us assume that the half commutator formula holds:
\begin{eqnarray}
 [Y(u,z_1), Y(v,z_2)]  &\simeq &
\mbox{Res}_{z_0} z_2^{-1} \delta \left(\frac{z_1 - z_0}
{z_2} \right) Y(Y(u,z_0)v,z_2) .\label{Hcommutator}
\end{eqnarray}
Then for $m\geq 0$ we have
\begin{eqnarray}
&&\mbox{Res}_{z_0} z_0^m z_2^{-1} \delta \left(\frac{z_1 - z_0}{z_2} \right)
Y(Y(u,z_0)v,z_2) \label{comm1}\\
&=&\mbox{Res}_{z_0} (z_1-z_2)^m\, z_2^{-1} \delta \left(\frac{z_1 -
z_0}{z_2} \right)
Y(Y(u,z_0)v,z_2) \label{comm2}\\
&=&(z_1-z_2)^m\, \mbox{Res}_{z_0} z_2^{-1} \delta \left(\frac{z_1 -
z_0}{z_2} \right)
Y(Y(u,z_0)v,z_2) \label{comm3}\\
&\simeq& (z_1-z_2)^m\,\left(Y(u,z_1) Y(v,z_2)
-\e_{u,v}\,Y(v,z_2)Y(u,z_1)\right) \label{comm4}\\
&=&(z_1-z_2)^m\,\mbox{Res}_{z_0}
z^{-1}_0\delta\left(\frac{z_1-z_2}{z_0}\right)
Y(u,z_1)Y(v,z_2)\nonumber\\
&&-\e _{u,v}\,(z_1-z_2)^m\,\mbox{Res}_{z_0}
z^{-1}_0\delta\left(\frac{z_2-z_1}{-z_0}\right)
Y(v,z_2)Y(u,z_1) \label{comm5}\\
&=&\mbox{Res}_{z_0} z_0^m
z^{-1}_0\delta\left(\frac{z_1-z_2}{z_0}\right)
Y(u,z_1)Y(v,z_2)\nonumber\\
&&-\e _{u,v}\,\mbox{Res}_{z_0} z_0^m
z^{-1}_0\delta\left(\frac{z_2-z_1}{-z_0}\right)
Y(v,z_2)Y(u,z_1), \label{comm6}
\end{eqnarray}
which is precisely the half Jacobi identity.

Here we use the usual properties of $\delta$-function. Since $m\geq 0$, we
have a
polynomial $z_0^m=((z_0-z_1)+z_1)^m$ which may be expanded in powers
$(z_0-z_1)^k z_1^{m-k}$. In the presence of $\delta$-function we may replace
a power $(z_0-z_1)^k$ by $(-z_2)^k$, and as a result we get a polynomial
$(z_1-z_2)^m$. In this way we get (\ref{comm2}) from (\ref{comm1}) and
(\ref{comm6}) from (\ref{comm5}).
Since the polynomial $(z_1-z_2)^m$ is a linear combination of powers
$z_1^kz_2^{m-k}$, $k\geq 0$, $m-k\geq 0$, to establish the equality of
principal parts of (\ref{comm3}) and (\ref{comm4}) it is enough to know
the equality of principal parts of formal Laurent series appearing
in the half commutator formula (\ref{Hcommutator}).\qed

\begin{lem}\label{assoc} In the definition of VLA-SS the half Jacobi
identity can be
equivalently substituted by the half associator formula.
\end{lem}

\pf Let us assume that the half associator formula holds:
\begin{eqnarray}
 Y(Y(u,z_0)v,z_2)  &\simeq &
\mbox{Res}_{z_1}z^{-1}_0\delta\left(\frac{z_1-z_2}{z_0}\right)
Y(u,z_1)Y(v,z_2)\nonumber\\
&&-\e _{u,v}\,\mbox{Res}_{z_1}
z^{-1}_0\delta\left(\frac{z_2-z_1}{-z_0}\right)
Y(v,z_2)Y(u,z_1).\label{Hassociator}
\end{eqnarray}
Let  $n\geq 0$. By using arguments as in the proof of Lemma \ref{comm} we
get
\begin{eqnarray}
&&\mbox{Res}_{z_1}z_1^n z_2^{-1}\delta\left(\frac{z_1-z_0}{z_2}\right)
Y(Y(u,z_0)v,z_2)\nonumber\\
&=&\mbox{Res}_{z_1}(z_2+z_0)^n\,
z_2^{-1}\delta\left(\frac{z_1-z_0}{z_2}\right)
Y(Y(u,z_0)v,z_2)\nonumber\\
&=&(z_2+z_0)^n\,Y(Y(u,z_0)v,z_2)\nonumber\\
&\simeq&
(z_2+z_0)^n\,\mbox{Res}_{z_1}z^{-1}_0\delta\left(\frac{z_1-z_2}{z_0}\right)
Y(u,z_1)Y(v,z_2)\nonumber\\
&&-\e _{u,v}\,(z_2+z_0)^n\,\mbox{Res}_{z_1}
z^{-1}_0\delta\left(\frac{z_2-z_1}{-z_0}\right)
Y(v,z_2)Y(u,z_1)\nonumber\\
&=&
\mbox{Res}_{z_1}(z_2+z_0)^n\,z^{-1}_0\delta\left(\frac{z_1-z_2}{z_0}\right)
Y(u,z_1)Y(v,z_2)\nonumber\\
&&-\e _{u,v}\,\mbox{Res}_{z_1}(z_2+z_0)^n\,
z^{-1}_0\delta\left(\frac{z_2-z_1}{-z_0}\right)
Y(v,z_2)Y(u,z_1)\nonumber\\
&=& \mbox{Res}_{z_1}z_1^n\,z^{-1}_0\delta\left(\frac{z_1-z_2}{z_0}\right)
Y(u,z_1)Y(v,z_2)\nonumber\\
&&-\e _{u,v}\,\mbox{Res}_{z_1}z_1^n\,
z^{-1}_0\delta\left(\frac{z_2-z_1}{-z_0}\right)
Y(v,z_2)Y(u,z_1),\nonumber
\end{eqnarray}
which is precisely the half Jacobi identity. \qed

\begin{rem}\label{dokazi2}
The proofs of Lemmas \ref{comm} and \ref{assoc} are modifications of the
proofs
of the corresponding statements for vertex superalgebras given in
[Li, Propositions 2.2.4 and 2.2.6].
\end{rem}

\begin{lem}\label{spancomm}
Let $U$ be a VLA-J-SS. Let \ $U_1$ be a subspace spanned by vectors of the
form
\begin{equation}\label{spancomm1}
u_mv_nw-\e_{u,v}\,v_nu_mw - \sum_{i \geq 0} {m \choose i}(u_iv)_{n+m-i}w \
\end{equation}
for all homogeneous $u,v,w\in U$ and $m, n\in \N$.
Let $U_2$ be a subspace of $U$ spanned by
vectors of the form
\begin{equation}\label{spanjacobi}
\begin{array}{c}
\displaystyle{\sum_{i\geq0} (-1)^i {k \choose i}(u_{m+k-i}(v_{n+i}w) -
\e_{u,v}\,(-1)^k v_{n+k-i}(u_{m+i}w)) }\\
\displaystyle{-\sum_{i \geq 0} {m \choose i}(u_{k+i}v)_{m+n-i}w}
\end{array}
\end{equation}
for all homogeneous $u,v,w\in U$ and $k, m, n\in \N$. Then $U_1=U_2$.
\end{lem}

\pf
For two Laurent series $A$ and $B$ let us write $A\equiv B$ if the
coefficients
of principal parts of $A$ and $B$ are equal modulo elements in $U_1$. By
definition
we have
\begin{eqnarray}
 [Y(u,z_1), Y(v,z_2)]w  &\equiv &
\mbox{Res}_{z_0} z_2^{-1} \delta \left(\frac{z_1 - z_0}
{z_2} \right) Y(Y(u,z_0)v,z_2)w .\nonumber
\end{eqnarray}
Then, by almost copying the proof of Lemma \ref{comm}, we get
\begin{eqnarray}
&&\mbox{Res}_{z_0} z_0^m z_2^{-1} \delta \left(\frac{z_1 - z_0}{z_2} \right)
Y(Y(u,z_0)v,z_2)w \nonumber\\
&\equiv &\mbox{Res}_{z_0} z_0^m
z^{-1}_0\delta\left(\frac{z_1-z_2}{z_0}\right)
Y(u,z_1)Y(v,z_2)w\nonumber\\
&&-\e _{u,v}\,\mbox{Res}_{z_0} z_0^m
z^{-1}_0\delta\left(\frac{z_2-z_1}{-z_0}\right)
Y(v,z_2)Y(u,z_1)w, \nonumber
\end{eqnarray}
which implies $U_2\subset U_1$. It is clear that $U_1\subset U_2$. \qed

If $U$ is a VLA-J-SS, we shall denote by $\langle\mbox{Jacobi\,}\rangle$ the
two-sided ideal in $U$ generated by vectors of the form (\ref{spanjacobi}),
or equivalently,
by vectors of the form (\ref{spancomm1}). It is clear that
$\langle\mbox{Jacobi\,}\rangle$ is $D$ invariant since $D$ is a derivation
for all multiplications. It is also clear that
$\langle\mbox{Jacobi\,}\rangle$ is a $\Z_2$-graded subspace of $U$ and that
we have:
\begin{lem}\label{U_J} Let $U$ be a VLA-J-SS and let
$$
U_J=U/\langle\mbox{Jacobi\,}\rangle.
$$
Then $U_J$ is a VLA-SS with the universal property that any
homomorphism $\varphi$ from $U$ to a vertex Lie superalgebra $V$ factors
through $U_J$
by an homomorphism $\varphi_{J}$ from $U_J$ to $V$.
\end{lem}

\begin{lem}\label{Wss} Let $W$ be a VLA-SS and let $\langle\mbox{skew
symm.}\rangle$
be a subspace spanned by vectors
\begin{equation}\label{spanss1}
v_nw+\e_{v,w}\,\sum_{k\geq 0}(-1)^{n+k}(D^k/k!)\,w_{n+k}v
\end{equation}
for all homogeneous $v,w\in W$ and $n\geq 0$. Then
$\langle\mbox{skew symm.}\rangle$ is a $\Z_2$-graded two-sided
ideal in $W$ invariant for $D$. Moreover,
$$
W_{ss}=W/\langle\mbox{skew symm.}\rangle
$$
is a VSLA with the universal property that any
homomorphism $\varphi$ from $W$ to a vertex Lie superalgebra $V$ factors
through $W_{ss}$
by an homomorphism $\varphi_{ss}$ from $W_{ss}$ to $V$.
\end{lem}

\pf  Set $W_1=\langle\mbox{skew symm.}\rangle$. It is clear that
$DW_1\subset W_1$
since by assumption $D$ is a derivation for all multiplications.

For two Laurent series $A$ and $B$ let us write $A\equiv B$ if the
coefficients
of principal parts of $A$ and $B$ are equal modulo elements in $W_1$. By
definition
we have
\begin{eqnarray}
&&e^{-z_2D}Y(v,z_2)w\equiv\e_{v,w}Y(w,-z_2)v.\nonumber
\end{eqnarray}
In order to show that $W_1$ is a left ideal consider
\begin{eqnarray}
&&Y(u,z_0)e^{-z_2D}Y(v,z_2)w\nonumber\\
&=&e^{-z_2D}Y(u,z_0+z_2)Y(v,z_2)w\nonumber\\
&=&e^{-z_2D}\mbox{Res}_{z_1} z_1^{-1} \delta \left(\frac{z_0 + z_2}
{z_1} \right) Y(u,z_0+z_2)Y(v,z_2)w\nonumber\\
&=&e^{-z_2D}\mbox{Res}_{z_1} z_0^{-1} \delta \left(\frac{z_1 - z_2}
{z_0} \right) Y(u,z_1)Y(v,z_2)w\nonumber\\
&\simeq& \e_{u,v}\,e^{-z_2D}\mbox{Res}_{z_1} z_0^{-1} \delta
\left(\frac{z_2 - z_1}
{-z_0} \right) Y(v,z_2)Y(u,z_1)w\nonumber\\
&&+e^{-z_2D}\mbox{Res}_{z_1} z_2^{-1} \delta \left(\frac{z_1 - z_0}
{z_2} \right) Y(Y(u,z_0)v,z_2)w\nonumber\\
&=& \e_{u,v}\,\mbox{Res}_{z_1} z_0^{-1} \delta \left(\frac{z_2 - z_1}
{-z_0} \right) e^{-z_2D}Y(v,z_2)Y(u,z_1)w\nonumber\\
&&+e^{-z_2D}Y(Y(u,z_0)v,z_2)w\nonumber\\
&\equiv& \e_{u,v}\,\e_{u,v}\,\e_{w,v}\,
\mbox{Res}_{z_1} (-z_2)^{-1} \delta \left(\frac{z_0 - z_1}
{-z_2} \right) Y(Y(u,z_1)w,-z_2)v\nonumber\\
&&+e^{-z_2D}Y(Y(u,z_0)v,z_2)w\nonumber\\
&\simeq& \e_{v,w}\,Y(u,z_0)Y(w,-z_2)v\nonumber\\
&&-\e_{u,w}\,\e_{v,w}\,Y(w,-z_2)Y(u,z_0)v\nonumber\\
&&+e^{-z_2D}Y(Y(u,z_0)v,z_2)w\nonumber\\
&\equiv& \e_{v,w}\,Y(u,z_0)Y(w,-z_2)v.\nonumber
\end{eqnarray}
Hence we have
\begin{eqnarray}
&&Y(u,z_0)\left(e^{-z_2D}Y(v,z_2)w-\e_{v,w}Y(w,-z_2)v\right)\equiv
0.\nonumber
\end{eqnarray}
This shows that $W_1$ is a left ideal, i.e. that $u_nW_1\subset W_1$ for all
$u\in W$
and $n\geq 0$. It remains to show that $W_1$ is a right ideal, so consider
\begin{eqnarray}
&&Y(Y(u,z_0)v,z_2)w\nonumber\\
&\equiv& \e_{u,w}\e_{v,w}e^{z_2D}Y(w,-z_2)Y(u,z_0)v\label{rideal1}\\
&\equiv&
\e_{u,w}\,\e_{v,w}\,e^{z_2D}Y(w,-z_2)\e_{v,u}\,e^{z_0D}Y(v,-z_0)u\label{ride
al2}\\
&\equiv& \e_{v,u}\,Y(e^{z_0D}Y(v,-z_0)u,z_2)w.\nonumber
\end{eqnarray}
Hence we have
\begin{eqnarray}
&&Y\left(Y(u,z_0)v-\e_{v,u}\,e^{z_0D}Y(v,-z_0)u,z_2\right)w\equiv
0.\nonumber
\end{eqnarray}
This shows that $W_1$ is a right ideal, i.e. that $(W_1)_nw\subset W_1$ for
all  $w\in W$
and $n\geq 0$. Note that we obtained (\ref{rideal2}) from (\ref{rideal1}) by
using the fact that $W_1$ is a left ideal invariant for $D$. The remaining
statements are clear. \qed

If $U$ is a VLA-J-SS, we set
$$
U_{J+ss}=\left( U_J\right)_{ss}=U/\langle\mbox{Jacobi\,}\rangle/
\langle\mbox{skew symm.}\rangle.
$$
Moreover, for  any homomorphism $\varphi$ from $U$ to a vertex Lie
superalgebra $V$
we set
$$
\varphi_{J+ss}=\left( \varphi_J\right)_{ss}=\varphi_{ss}\circ\varphi_J.
$$
Then clearly $\varphi : U\to V$ factors through $U_{J+ss}$ by the
homomorphism
$\varphi_{J+ss} : U_{J+ss}\to V$.
Obviously
$$
U_{J+ss}=U/\langle J+ss\rangle ,
$$
where $\langle J+ss\rangle$ is a two-sided ideal in $U$ generated by
elements
of the form (\ref{spancomm1}) and (\ref{spanss1}). It is clear that
$\langle J+ss\rangle$ is $D$ invariant since $D$ is a derivation for all
multiplications. It is also clear that $\langle J+ss\rangle$ is
$\Z _2$-graded.


\section{Vertex Lie algebras generated by formulas}

Let $S$ be a vector space given bilinear maps $F_n^k : S\times S\to S$ for
$n,k\geq 0$.
Than we can ask under what conditions on $F_n^k$ there is a vertex algebra
$V$ containing
$S$ and generated by the set of fields $\{Y(u,z)\mid u\in S\}$ for which the
commutator
formula is
$$
 [Y(u,z_1), Y(v,z_2)]
       = \sum_{n,k \geq 0} \frac{(-1)^n}{n!} \left(\frac{d}{dz_1} \right)^n
z^{-1}_2
\delta(z_1/z_2)\left(\frac{d}{dz_2} \right)^kY( F_n^k(u,v),z_2) .
$$
If there is such $V$, we could say that $V$ is generated by a formula
defined
by the set of maps $\{F_n^k\mid n,k\geq 0\}$.
There are many examples when this is the case; besides the examples in [Li]
mentioned
before one may consult [DLM] and the references therein.

We shall say that a $\Z_2$-graded vector space $S$ is a formula
if it is equipped with infinitely many $\Z_2$-grading preserving linear maps
$$
F_n : S\otimes S\to \C[D]\otimes S,\ \ \ \ n\in \N
$$
with the space $\C[D]$ of formal polynomials in $D$ considered to be even,
such that for  $u,v\in S$
\begin{equation}\label{nvelik0}
F_n(u,v)=0\ \ \
\mbox{for}\ \ n \ \ \mbox{sufficiently\ large}.
\end{equation}
We shall write
$$
F_n(u,v)=u_nv=\sum_{k\geq 0}D^k\otimes F_n^k(u,v),\ \ \ u,v, F_n^k(u,v)\in
S.
$$
It is clear that giving the maps $F_n$ for all $n\geq 0$ is equivalent to
giving
a linear map
$$
Y : S\otimes S\to z^{-1}(\C[D]\otimes S)[[z^{-1}]], \ \ \
Y(u,z)v=\sum_{n\geq 0}\, F_n(u,v)z^{-n-1}
$$
and it would be proper to say that a formula is a given map $Y_S=Y$, or
a pair $(S,Y)$, with the prescribed
properties.

We define a linear operator $D(D^k\otimes u)=D^{k+1}\otimes u$
and consider $ S\subset \C[D]\otimes S$. Then $D^k\otimes u$ can also be
written as
$D^k u$, and in general elements in $\C[D]\otimes S$ can
be written as linear combinations of elements of the form $P(D)u$, where
$u\in U$ and
$P(D)$ is a polynomial $P$ of the operator $D$. We shall use the
notation $A|_{D=0}=0$ with the obvious meaning that $A\in D\C[D]\otimes S$.

{\bf Example 1.}\ \ Let $S=S^0={\bf g} \oplus \C c$ be a sum of
1-dimensional space and a Lie algebra
{\bf g} with an invariant bilinear form $\langle\cdot,\cdot \rangle$. Set
\begin{eqnarray*}
&& x_0y=[x,y], \ \ \ x_0c=c_0x=c_0c=0,\\
&& x_1y=\langle x,y\rangle c, \ \ \ x_1c=c_1x=c_1c=0
\end{eqnarray*}
for all $x,y\in {\bf g}$. Set $F_n=0$ for $n\geq 2$. Then $S$ is a formula.
We could also write

$$
Y(x,z)y=\frac{[x,y]}{z}+\frac{\langle x,y\rangle c}{z^2}, \ \ \ \
Y(x,z)c=0, \ \ \ \
Y(c,z)=0.
$$

{\bf Example 2.}\ \ Let $S=S^0=\C\omega \oplus \C c$ be a
2-dimensional space. Set
$$
Y(\omega,z)\omega =\frac{D\omega}{z}+\frac{2\omega}{z^2}+\frac{(1/2)c}{z^4},
\ \ \ \
Y(\omega,z)c=0, \ \ \ \
Y(c,z)=0.
$$
Then $S$ is a formula.

\begin{lem}\label{prosirenjeop}
Bilinear operations $F_n$ on $S\subset \C[D]\otimes S$ extend in a unique
way to
$\Z_2$-grading preserving linear maps
$$
F_n : \left(\C[D]\otimes S\right)\otimes \left(\C[D]\otimes S\right)
\to \C[D]\otimes S,\ \ \ \ n\in \N
$$
such that
\begin{eqnarray}
&&(DA)_nB=-nA_{n-1}B ,\label{formder1}\\
&&D(A_nB)=(DA)_{n}B+A_n(DB) \label{formder2}
\end{eqnarray}
for all $A,B\in \C[D]\otimes S$ and $n\geq 0$, where we write
$F_n(A,B)=A_nB$.
Moreover, for $A,B\in \C[D]\otimes S$
$$
A_nB=0\ \ \
\mbox{for}\ \ n \ \ \mbox{sufficiently\ large}.
$$
\end{lem}

\pf First note that our conditions can be written as
$$
Y(DA,z)B=(d/dz)Y(A,z)B,\ \ \ \ Y(A,z)DB=(D-d/dz)Y(A,z)B.
$$
Hence it must be
\begin{equation}\label{formzaY}
Y(P(D)u,z)Q(D)v=P(d/dz)Q(D-d/dz)Y(u,z)v
\end{equation}
for $u,v\in S$ and any polynomials $P$ and $Q$. It is easy to check all
the properties stated. \qed

In a way Lemma \ref{prosirenjeop} means that a formula extends to a
VLA-J-SS.
If $\varphi : S\to V$ is a linear map from a formula $S$ to a vertex Lie
superalgebra $V$, we extend it to a map
$$
\tilde\varphi : \C[D]\otimes S\to V,\ \ \ \
\tilde\varphi(D^k\otimes u)=D^k \varphi(u).
$$
We shall say that $\varphi : S\to V$ is a homomorphism from $S$ to $V$ if
$\tilde\varphi$ is a homomorphism (cf. Section 6). By using (\ref{formzaY})
it is
easy to see that a linear map $\varphi : S\to V$ is a homomorphism if and
only if
$$
\tilde\varphi(u_nv)=\varphi(u)_n\varphi(v)\ \ \ \ \ \mbox{for all}\ \
u,v\in S.
$$
Obviously Lemmas \ref{U_J} and \ref{Wss} imply:
\begin{lem}\label{homoS} Any homomorphism $\varphi : S\to V$ from a formula
$S$
to a vertex Lie superalgebra $V$
factors through $(\C[D]\otimes S)_{J+ss}$ by a homomorphism
$\tilde\varphi_{J+ss} : (\C[D]\otimes S)_{J+ss} \to V$.
\end{lem}

We shall say that $S$ is an $(J+ss)$-injective formula if the restriction of
the quotient map to $S$:
$$
S\hookrightarrow \C[D]\otimes S\to (\C[D]\otimes S)_{J+ss}
$$
is an injection.
\begin{lem}\label{karakinj}
A formula $S$ is $(J+ss)$-injective if and only if there is a vertex Lie
superalgebra $V$
and an injective homomorphism $\varphi$ from $S$ to $V$.
\end{lem}
\pf
If $S$ is $(J+ss)$-injective, then take $V=(\C[D]\otimes S)_{J+ss}$.
If $\varphi : S\to V$ is injective, then the restriction on $S$ of the
quotient map
must be injective. \qed

In the case when $S$ is $(J+ss)$-injective we can
consider $S$ as a subspace of the vertex Lie superalgebra\, $U=(\C[D]\otimes
S)_{J+ss}$
for which the multiplications $u_nv$ are determined by the
formula
\begin{equation}\label{formtm1}
u_nv=\sum_{k\geq 0}D^k F_n^k(u,v),\ \ \ u,v, F_n^k(u,v)\in S.
\end{equation}
Moreover, for a VLSA $V$ a linear map $\varphi : U_{J+ss}\to V$
is a homomorphism if and only if
$$
\varphi(u)_n\varphi(v)=
\sum_{k\geq 0}D^k \varphi\left(F_n^k(u,v)\right)
\quad\mbox{for}\quad u,v\in S.
$$

In another words, relations (\ref{formtm1}) consistently determine the
multiplications $u_nv$ in a vertex Lie superalgebra
$U=\mbox{span}\,\{D^kw\mid w\in S, k\geq 0\}$ if and only if $S$ is
$(J+ss)$-injective.
What remains to be seen is under what conditions this is the case. We have
only a partial answer to this problem which, at least, covers the
examples listed before.

The following two lemmas are the consequence of (\ref{formzaY}), some
calculations get to be simpler by using [Li, (2.1.8) and (2.1.15)]:
\begin{lem}\label{formidealcomm}
For $u,v,w\in S$ and polynomials $P$, $Q$ and $R$ we have in  $\C[D]\otimes
S$
\begin{eqnarray*}
&&[Y(P(D)u,z_1), Y(Q(D)v,z_2)]R(D)w \\
&&-\mbox{Res}_{z_0} z_2^{-1} \delta \left(\frac{z_1 - z_0}
{z_2} \right) Y(Y(P(D)u,z_0)Q(D)v,z_2)R(D)w \\
&=&P(d/dz_1)Q(d/dz_2)R(D-d/dz_1-d/dz_2)\\
&&\left([Y(u,z_1), Y(v,z_2)]w
-\mbox{Res}_{z_0} z_2^{-1} \delta \left(\frac{z_1 - z_0}
{z_2} \right) Y(Y(u,z_0)v,z_2)w\right).
\end{eqnarray*}
\end{lem}
\begin{lem}\label{formidealss}
For $u,v\in S$ and polynomials $P$ and $Q$ we have in  $\C[D]\otimes S$
\begin{eqnarray*}
&&Y(P(D)u,z)Q(D)v-\e_{u,v}e^{zD}Y(Q(D)v,-z)P(D)u\\
&=&P(d/dz)Q(D-d/dz)\left(Y(u,z)v-\e_{u,v}e^{zD}Y(v,-z)u\right)\\
&\equiv& P(d/dz)Q(-d/dz)\left(Y(u,z)v-\e_{u,v}Y(v,-z)u\right)|_{D=0},
\end{eqnarray*}
where $\equiv$ means that the coefficients in the principal part of Laurent
series
in $z$ are equal modulo subspace $D\C[D]\otimes S$.
\end{lem}
\begin{prop}\label{Liealg}
Let $S$ be a formula and assume that $F_n(u,v)\in S$
for all $u,v\in S$ and $n\geq 0$.
Then $\C[D]\otimes S$ is a vertex Lie superalgebra if and only if
$F_n=0$ for all $n\geq 1$ and $F_0$ defines a Lie superalgebra structure on
$S$.
\end{prop}
\pf
Let  $\C[D]\otimes S$ be a VLSA. Then the assumption $F_n(u,v)\in S$
together
with the half skew symmetry relation (\ref{HSScomp}) implies $F_n=0$ for all
$n\geq 1$
and $u_0v=-\e_{u,v}\,v_0u$. The Jacobi identity (\ref{jacobicomp}) for
$k=m=n=0$
is the Jacobi identity for Lie superalgebras. In the same way we see that
the
assumption $F_n=0$ for all $n\geq 1$ and $F_0$ a Lie superalgebra commutator
implies the half Jacobi identity and the half skew symmetry relation for
elements in $S$, so then by Lemmas \ref{formidealcomm} and \ref{formidealss}
they hold for all elements in $\C[D]\otimes S$.  \qed

Clearly Lemmas \ref{comm} and \ref{formidealcomm} imply the following:
\begin{prop}\label{Jformula}
Let $S$ be a formula. Then $\C[D]\otimes S$ is a VLA-SS if and only if the
half
commutator relation (\ref{Hcommutator}) holds in $\C[D]\otimes S$
for all $u,v,w \in S$.
\end{prop}

Note that by Lemmas \ref{spancomm} and \ref{formidealcomm} the ideal
$\langle\mbox{Jacobi\,}\rangle$ in $\C[D]\otimes S$ is a two-sided $D$
invariant ideal generated by elements
\begin{equation}\label{spancommform}
u_mv_nw-\e_{u,v}\,v_nu_mw - \sum_{i \geq 0} {m \choose i}(u_iv)_{n+m-i}w \
\end{equation}
for all homogeneous $u,v,w\in S$ and $m, n\in \N$. Also note that
by Lemmas \ref{Wss} and \ref{formidealss} the ideal
$\langle\mbox{skew symm.}\rangle$ in $(\C[D]\otimes S)_J$ is a two-sided $D$
invariant ideal spanned by elements
\begin{equation}\label{spanformss1}
u_nv+\e_{u,v}\,\sum_{k\geq 0}(-1)^{n+k}(D^k/k!)\,v_{n+k}u
\end{equation}
for all homogeneous $u,v\in S$ and $ n\in \N$. Hence it is clear from
our construction of
$$
(\C[D]\otimes S)_{J+ss}=(\C[D]\otimes S)/\langle\mbox{Jacobi\,}\rangle/
\langle\mbox{skew symm.}\rangle=(\C[D]\otimes S)/\langle J+ss\rangle
$$
that we have:
\begin{prop}\label{forminj}
Let $S$  be a formula such that $\langle J+ss\rangle\subset D\C[D]\otimes
S$,
or equivalently, such that
\begin{equation}\label{forminjJ}
\langle\mbox{Jacobi\,}\rangle\subset D\C[D]\otimes S
\end{equation}
and that for all homogeneous $u,v\in S$ and $n\geq 0$
\begin{equation}\label{forminjss}
u_nv+\e_{u,v}(-1)^n\,v_nu\in  D\C[D]\otimes S.
\end{equation}
Then $S$ is $(J+ss)$-injective.
\end{prop}

Note that the condition (\ref{forminjss}) is equivalent to
$F_n^0(u,v)=-\e_{u,v}(-1)^n\,F_n^0(v,u)$.
In Example 1 it reads $[x,y]=-[y,x]$,
$\langle x,y\rangle=\langle y,x\rangle$, whereas in Example 2 it is
obviously
satisfied.

In Example 1 the condition (\ref{forminjJ}) holds since
$\langle\mbox{Jacobi\,}\rangle=0$.
This amounts to a verification of relations
\begin{eqnarray*}
&&u_0v_0w-v_0u_0w = (u_0v)_{0}w, \\
&&u_0v_1w-v_1u_0w = (u_0v)_{1}w, \\
&&u_1v_0w-v_0u_1w = (u_0v)_{1}w +(u_1v)_{0}w
\end{eqnarray*}
for $u,v,w\in S$, the relation $u_1v_1w-v_1u_1w =  (u_1v)_{1}w$ and all the
other
obviously hold.

In Example 2 the condition (\ref{forminjJ}) holds
because $\langle\mbox{Jacobi\,}\rangle=D\C[D]c$. In principle we could
see this by a direct computation of elements in (\ref{spancommform}),
i.e. by a direct computation of the principal part of
$$
[Y(\omega,z_1), Y(\omega,z_2)]\omega
-\sum_{i=0}^1 \frac{(-1)^i}{i!} \left(\frac{d}{dz_1} \right)^i z^{-1}_2
\delta(z_1/z_2)Y(\omega_i\omega,z_2)\omega,
$$
but since we know there is a Virasoro vertex operator algebra $V$ of level
$\ell\neq 0$,
we can argue the other way around: First
note that $D\C[D]c$ is a two-sided $D$ invariant ideal in $\C[D]\otimes S$.
By a direct computation in $\C[D]\otimes S$ we see that
$$
\omega_0\omega_3\omega-\omega_3\omega_0\omega-
(\omega_0\omega)_3\omega
=-(1/2) Dc,
$$
so $\langle\mbox{Jacobi\,}\rangle\supset D\C[D]c$.
A linear map $\C[D]\otimes S\to V$
defined by $c\mapsto \ell\,\1$, $\omega\mapsto \omega$ factors through the
quotient
$$
\C c\oplus \C[D]\omega\cong (\C[D]\otimes S)/ D\C[D]c\to
\mbox{span}\,\{\1, D^k\omega\mid k\geq 0\}\subset V
$$
and it is obviously a homomorphism of VLA-J-SS. Since $\{\1, D^k\omega\mid
k\geq 0\}$
is linearly independent set in $V$, this map is an isomorphism, and hence
$\langle\mbox{Jacobi\,}\rangle=D\C[D]c=D\C[D]\otimes c$.

{\bf Example 3.}\ \
Let $S=S^0=B \oplus \C c$ be a sum of 1-dimensional space and an algebra
$(B,\cdot)$ with a symmetric bilinear form $\langle\cdot,\cdot\rangle$.
For $u,v\in B$ set
$$
Y(u,z)v=\frac{D(u\cdot v)}{z}+\frac{u\cdot v+v\cdot u}{z^2}+
\frac{\frac{1}{2}\langle u,v\rangle c}{z^4}, \qquad
Y(u,z)c=0, \ Y(c,z)=0.
$$
It is easy to see that for $u,v,w\in S$ the elements of the form
(\ref{spancommform})
and (\ref{spanformss1}) are in the ideal $D\C[D]\otimes c$,
or equivalently, that $\langle J+ss\rangle$ is contained in the ideal
$D\C[D]\otimes c$, if and only if $B$ is a right Novikov algebra (cf. [O]),
i.e. if and only if
$$
u\cdot(v\cdot w)=v\cdot(u\cdot w),
$$
$$
(v\cdot w)\cdot u+v\cdot(u\cdot w)=v\cdot(w\cdot u)+(v\cdot u)\cdot w,
$$
$$
\langle u\cdot v,w\rangle =\langle v\cdot u,w\rangle
=\langle v,u\cdot w\rangle =\langle v,w\cdot u\rangle.
$$
If $B$ is an associative commutative algebra with a symmetric associative
bilinear form $\langle\cdot,\cdot\rangle$, then all these conditions are
clearly
satisfied. For an obvious noncommutative example we take a linear functional
$\lambda$ on a vector space $B$ and set $u\cdot v=\lambda (u)v$,
$\langle u, v\rangle=\lambda(u)\lambda(v)$. Note that $\lambda(\omega)=1$
implies $\omega\cdot\omega=\omega$.

If the above conditions are satisfied, then $S$ is (J+ss)-injective, that is
$S\subset U=(\C[D]\otimes S)_{J+ss}$. If we set
$Y(u,z)=\sum_{n\in\Z}u(n)z^{-n-2}$,
then in the corresponding Lie algebra $\L(U)$ we have the commutation
relations
$$
[u(n),v(m)]=(n+1)(v\cdot u)(n+m)-(m+1)(u\cdot v)(n+m)
+\frac12{n+1 \choose 3}\delta_{n+m,0}\langle u, v\rangle c.
$$

We can summarize our constructions in the following:
\begin{thm}\label{formtmZ}
Let $(S,Y_S)$ be a formula given by
$$
Y_S(u,z)v=\sum_{n,k\geq 0}\frac{D^k\otimes F_n^k(u,v)}{z^{n+1}}
$$
and assume  that
$\langle\mbox{Jacobi\,}\rangle|_{D=0}\,=0$ and
$F_n^0(u,v)=-\e_{u,v}(-1)^n\,F_n^0(v,u)$
for $u,v\in S$ and $n\geq 0$.

Then there exists a vertex superalgebra $V$ and an injective homomorphism
$ S\hookrightarrow V$ such that $V$ is generated by the set of fields
$\{Y(u,z)\mid u\in S\}$ for which the commutator formula is
\begin{equation}\label{11}
 [Y(u,z_1), Y(v,z_2)]
       = \sum_{n,k \geq 0} \frac{(-1)^n}{n!} \left(\frac{d}{dz_1} \right)^n
z^{-1}_2
\delta(z_1/z_2)\left(\frac{d}{dz_2} \right)^kY( F_n^k(u,v),z_2) .
\end{equation}
Among all such vertex superalgebras, i.e. $ S\hookrightarrow V$ to be
precise,
there exists a vertex superalgebra
$\V(\langle S\rangle)$ such that any other $V$ is a quotient of
$\V(\langle S\rangle)$.
Moreover, for any $V$ as above we have
\begin{equation}\label{12}
Y_V(u,z)v\simeq \sum_{n,k\geq 0}\frac{D^k F_n^k(u,v)}{z^{n+1}}.
\end{equation}
\end{thm}

\pf By Proposition \ref{forminj} the formula $S$ is $(J+ss)$-injective.
So consider $S\subset\langle S\rangle$,
where $\langle S\rangle=(\C[D]\otimes S)_{J+ss}$ is a vertex Lie algebra.
Note that by our construction on the quotient
$\langle S\rangle$ we have a relation
\begin{equation}\label{13}
u_nv=\sum_{k\geq 0}D^k F_n^k(u,v)
\end{equation}
for $u,v\in S$, with $F_n^k(u,v)\in S$. Let $\V(\langle S\rangle)$ be the
universal enveloping vertex algebra of $\langle S\rangle$ and, as usual,
consider $S\subset \langle S\rangle\subset \V(\langle S\rangle)$. Since
$\langle S\rangle$ is generated by $S$ and $D$, vertex superalgebra
$\V(\langle S\rangle)$ is generated by $\{Y(u,z)\mid u\in S\}$ and
(\ref{11})
is just the commutator formula in which the product $u_nv$ in
$\langle S\rangle\subset \V(\langle S\rangle)$ is expressed by using
(\ref{13}).
Since on $\langle S\rangle\subset \V(\langle S\rangle)$ the operations
for $n\geq 0$ coincide, (\ref{12}) follows as well.

Now assume that $V$ is another vertex superalgebra and that
$\varphi : S\hookrightarrow V$ is an injective
homomorphism from $S$ to $V$. By Lemma \ref{homoS} the homomorphism
$\varphi$
factors through $\langle S\rangle$ by an homomorphism
$\tilde\varphi_{J+ss} : \langle S\rangle\to V$. By the universal property of
$\V(\langle S\rangle)$ this homomorphism $\tilde\varphi_{J+ss}$ extends to a
homomorphism $\V(\langle S\rangle)\to V$, and it is surjective since $S$
generates $V$. \qed

We may illustrate Theorem \ref{formtmZ} on the Virasoro formula $S$ given in
Example~2.
We saw that $U=(\C[D]\otimes S)_{J+ss}$ is generated by $S$ and that
$\langle S\rangle =U=\C c\oplus \C[D]\omega$ with $Dc=0$. The corresponding
Lie algebra
$\L(U)$ is the Virasoro Lie algebra
$$
\L(U)=\C c\oplus\mbox{span}\,\{L(n)\mid n\in\Z\}
$$
with a central element $c$ and the commutation relations
$$
[L(n),L(m)]=(n-m)L(n+m)+\frac12{n+1 \choose 3}\delta_{n+m,0}c.
$$
Here we use the usual notation
$$
Y(\omega,z)=\sum_{n\in\Z} \omega_nz^{-n-1}=\sum_{n\in\Z} L(n)z^{-n-2},
$$
i.e. $L(n)= \omega_{n+1}$, and $c=c_{-1}$, $c_n=0$ for $n\neq-1$. We also
have
$$
\L_-(U)=\C c\oplus\mbox{span}\,\{L(n)\mid n\leq -2\},\qquad
\L_+(U)=\mbox{span}\,\{ L(n)\mid n\geq -1\}.
$$
The derivation $D$ of the Lie algebra $\L(U)$, defined by
$D(u_n)=(Du)_n$, equals $\mbox{ad}\,(L(-1))$, i.e.
$$
D(L(n))=-(n+1)L(n-1)=[L(-1),L(n)],\qquad D(c)=0.
$$
By construction, the universal enveloping vertex algebra $\V(\langle
S\rangle)$
of the vertex Lie algebra $U=\langle S\rangle$ can be viewed as the
universal enveloping algebra $\U(\L_-(U))$ of the Lie algebra $\L_-(U)$.
Under this identification $\1$ is the identity and $D$ is a derivation
of the associative algebra $\U(\L_-(U))$ extending the derivation $D$
of the Lie algebra $\L_-(U)$. The vertex algebra $\V(\langle S\rangle)$ is
generated
by fields
$$
Y(\omega,z), \quad Y(c,z)=c\quad\mbox{and} \quad
Y(\1,z)=\mbox{id}\,_{\V(\langle S\rangle)}.
$$
Note that $c$ acts as a left multiplication by $c$, and that $D=L(-1)$.

For any vertex operator superalgebra $V$ of level $\ell$ with a
conformal vector $\omega$ we have an
injective homomorphism $ S\hookrightarrow V$, defined by $c\mapsto
\ell\,\1$,
$\omega\mapsto \omega$, which extends to a homomorphism
$ \V(\langle S\rangle)\rightarrow V$ of vertex superalgebras. Clearly
$\langle c-\ell\,\1\rangle=(c-\ell\,\mbox{id})\V(\langle S\rangle)$ is an
ideal
in $\V(\langle S\rangle)$, so we have a  vertex superalgebra homomorphism
$$
\V(\langle S\rangle)/\langle c-\ell\,\1\rangle\rightarrow V.
$$


\section{Conformal vectors}

In the previous section we have considered three examples, and in each of
them the corresponding
vertex Lie algebra $U=(\C[D]\otimes S)_{J+ss}$ is a quotient of
$\C[D]\otimes S$ by the
ideal $\langle J+ss\rangle=D\C[D]\otimes c$. In this section we will show
that this holds in general for certain class of formulas, roughly speaking
the ones corresponding to vertex operator superalgebras. For this reason we
need the notions of graded formulas and conformal vectors.

Let $(S,Y)$ be a formula and $d : S\to S$ an even linear map such that for
$i\in \Bbb Z_2$
$$
S^i=\oplus_{\lambda\in \Bbb R_+} S^i_\lambda\,, \qquad
dv=\lambda v \quad \mbox{for} \quad v\in S^i_\lambda .
$$
We extend $d$ to an even linear map on $\Bbb C[D]\otimes S$ by setting
$$
d\,(D^kv)=(\lambda+k)v \quad \mbox{for} \quad v\in S_\lambda, \ k\geq 0 .
$$
Here we write $S_\lambda=S^0_\lambda+S^1_\lambda$. So the space
$\C[D]\otimes S$
is graded by eigenspaces of $d$:
$$
\C[D]\otimes S=  \oplus_{\mu\in \Bbb R_+} (\C[D]\otimes S)_\mu
$$
\begin{equation}\label{821}
(\C[D]\otimes S)_\mu=\oplus_{k\geq 0,\,\lambda\geq 0,\,k+\lambda=\mu}\,
D^k\otimes S_\lambda.
\end{equation}
We shall say that the formula $(S,Y)$ is graded by $d$ \, if
\begin{equation}\label{83}
[d,Y(u,z)]v=(z\frac{d}{dz}+\lambda)Y(u,z)v \quad \mbox{for}
\quad u\in S_\lambda, \ v\in S.
\end{equation}
\begin{lem}\label{lem81}
For $A\in (\Bbb C[D]\otimes S)_\lambda$ and $B\in (\Bbb C[D]\otimes S)_\mu$
we have
\begin{equation}\label{84}
A_nB\in (\Bbb C[D]\otimes S)_{\lambda+\mu-n-1}.
\end{equation}
\end{lem}
\pf
By assumption (\ref{83}) the statement (\ref{84}) holds for $A=u$ and $B=v$
in $S$.
Since
$$
(DA)_nB=-nA_{n-1}B
$$
and
\begin{equation}\label{85}
D(\Bbb C[D]\otimes S)_{\lambda}\subset (\Bbb C[D]\otimes S)_{\lambda +1},
\end{equation}
we see by induction that  (\ref{84}) holds for $A\in \Bbb C[D]\otimes S$ and
$B=v\in S$.
Now the lemma follows by induction for all $B$ by using (\ref{85}) and
(\ref{formder2}),
i.e.
$$
A_n(DB)=(DA)_{n}B-D(A_nB). \ \ \ \Box
$$
Note that (\ref{85}) is equivalent to
\begin{equation}\label{86}
[d,D]=D.
\end{equation}

Let $(S,Y)$ be a formula and let $\omega$ and $c$ be two even nonzero
elements in $S$.
We shall say that $\omega$ is a conformal vector in $S$ with a central
element $c$ if
$$
Y(\omega,z)\omega =\frac{D\omega}{z}+\frac{2\omega}{z^2}+\frac{(1/2)c}{z^4},
$$
\begin{equation}\label{87}
Y(c,z)v=0, \quad Y(v,z)c=0 \qquad \mbox{for all}\quad v\in S,
\end{equation}
and if the formula $S$ is graded by a map $d$\, such that
$$
S_0=\Bbb C c,\qquad \mbox{dim}\, S_\lambda <\infty \quad\mbox{for}\quad
\lambda\in\Bbb R_+,
$$
\begin{equation}\label{89}
Y(\omega,z)v =\frac{Dv}{z}+\frac{dv}{z^2}+\frac{0}{z^3}+\dots
\quad\mbox{for}\quad v\in S.
\end{equation}
Clearly (\ref{87}) implies the following:
\begin{lem}\label{lem2} $D\Bbb C[D]\otimes c$ is a two-sided
$D$ invariant ideal in $\Bbb C[D]\otimes S$.
\end{lem}
\begin{lem}\label{lem3}
On $\Bbb C[D]\otimes S$ we have $\omega_0=D$ and $\omega_1=d$.
\end{lem}
\pf
Let $v\in S$. By the assumption (\ref{89}) we have
$\omega_0v=Dv$ and $\omega_1v=dv$, and for $u=D^kv$, $k\geq 1$, the relation
(\ref{formzaY}) implies
\begin{eqnarray*}
Y(\omega,z)D^kv &=& \left(D-\frac{d}{dz}\right)^k
\left(\frac{Dv}{z}+\frac{dv}{z^2}+\dots\right)\\
&=&\frac{D^{k+1}v}{z}+\frac{D^kdv+kD^kv}{z^2}+\dots .
\end{eqnarray*}
Hence \, $\omega_0 D^kv=DD^{k}v$, and (\ref{86}) implies \,
$\omega_1D^kv=dD^kv$.\qed
\begin{thm}\label{thm1}
Let $(S,Y)$ be a formula and let $\omega$ be a conformal vector in $S$
with a central element $c$. Then the following three conditions are
equivalent:

(i) $S$ is (J+ss)-injective.

(ii) $\langle J+ss\rangle = D\Bbb C[D]\otimes c$.

(iii) The elements of the form
\begin{equation}\label{881}
u_mv_nw-\e_{u,v}\,v_nu_mw - \sum_{i \geq 0} {m \choose i}(u_iv)_{n+m-i}w, \
\end{equation}
\begin{equation}\label{882}
v_nw+\e_{v,w}\,\sum_{k\geq 0}(-1)^{n+k}(D^k/k!)\,w_{n+k}v ,
\end{equation}

are in $D\Bbb C[D]\otimes c$ for all homogeneous $u,v,w\in S$ and $m, n\in
\N$.
\end{thm}
\pf
(i)$\Rightarrow$(ii).
Let $S$ be (J+ss)-injective and let
$U=(\Bbb C[D]\otimes S)_{J+ss}=\sum_{k\geq 0} D^kS$. As usual set
$Y_U(\omega,z)=\sum_{n\geq 0}\omega_nz^{-n-1}$.
Then for $v\in S_\lambda\subset U$ and $k\geq 1$ our assumptions imply
$\omega_0=D$ and $\omega_1(D^kv)=(\lambda+k)D^kv$.

The commutator formula for $\omega_0, \omega_1, \omega_2$ shows that on $U$
we
have a representation of the Lie algebra $sl_2$, in standard notation for
basis elements
$$
e=-\omega_2,\quad h=-2\omega_1,\quad f=\omega_0=D.
$$
By our assumptions $v\in S_\lambda\backslash\{0\}$ is a highest weight
vector
of $h$-weight $-2\lambda$, so $\Bbb C[D]v_\lambda$ is an irreducible Verma
module for
$\lambda>0$. Hence
$$
U=\Bbb C[D]c+\left(\oplus_{k\geq 0,\,\lambda>0} D^kS_\lambda\right),
$$
\begin{equation}\label{8xx}
\mbox{dim}\,(D^kS_\lambda)=\mbox{dim}\,S_\lambda\quad\mbox{for}\quad
\lambda>0.
\end{equation}
Since in $\C[D]\otimes S$ we have
$
\omega_0\omega_3\omega-\omega_3\omega_0\omega-
(\omega_0\omega)_3\omega
=-(1/2) Dc\in \langle\mbox{Jacobi\,}\rangle
$,
on the quotient $U$ we have $Dc=0$ and hence
$$
U=\Bbb Cc\oplus\left(\oplus_{k\geq 0,\,\lambda>0} D^kS_\lambda\right).
$$
In particular, $U=\oplus_{\mu\geq 0} U_\mu$ is graded
by eigenspaces of $\omega_1$:
\begin{equation}\label{883}
U_0=\Bbb Cc,\qquad
U_\mu=\oplus_{k\geq 0,\,\lambda>0,\,k+\lambda=\mu} D^kS_\lambda
\quad\mbox{for}\quad \mu>0.
\end{equation}
Since the two-sided ideal $\langle J+ss\rangle$ in $\C[D]\otimes S$
is generated by elements of the form (\ref{881}) and (\ref{882}), Lemma
\ref{lem81}
implies that $\langle J+ss\rangle$ is graded, i.e. invariant for
$\omega_1=d$,
and hence
\begin{equation}\label{884}
U_\mu\cong (\C[D]\otimes S)_\mu/\langle J+ss\rangle_\mu.
\end{equation}
It follows from (\ref{821}), (\ref{8xx}) and (\ref{883}) that
$$
\mbox{dim}\,U_\mu=\mbox{dim}\,(\C[D]\otimes S)_\mu-1 \quad\mbox{for}\quad
\mu>0,\,\mu\in\N,
$$
$$
\mbox{dim}\,U_\mu=\mbox{dim}\,(\C[D]\otimes S)_\mu \quad\mbox{for}\quad
\mu>0,\,\mu\notin\N,
$$
so (\ref{884}) implies $\mbox{dim}\,\langle J+ss\rangle_\mu= 1$ for $\mu>0$,
$\mu\in\N$, and
$\mbox{dim}\,\langle J+ss\rangle_\mu= 0$ for $\mu>0$, $\mu\notin\N$.
Hence  $D\C[D]\otimes c\subset\langle J+ss\rangle$ implies that
$$
\langle J+ss\rangle=D\C[D]\otimes c.
$$

(ii)$\Rightarrow$(iii) is obvious since by definition $\langle J+ss\rangle$
is generated by elements of the form (\ref{881}) and (\ref{882}).

(iii)$\Rightarrow$(i). Since $D\C[D]\otimes c$ is a two-sided ideal in
$\C[D]\otimes S$,
the assumption (iii) clearly implies $\langle J+ss\rangle\subset
D\C[D]\otimes c$.
Hence (i) follows by Proposition \ref{forminj}. \qed

\begin{rem}\label{confform}
As a conclusion we could say that in the case of (J+ss)-injective formulas
$S$
with a conformal vector $\omega$ it is enough to consider a VLA-J-SS \ $U$
of the
form
$$
U=\C c\oplus \left(\C[D]\otimes S'\right),\qquad S'=\oplus_{\lambda>0}\,
S_\lambda,
$$
defined by a map $Y$ on \ $S'\times S'$
$$
Y(u,z)v=\sum_{n\geq 0}\frac{u_nv}{z^{n+1}},\qquad u,v\in S',\quad u_nv\in U,
$$
$Y$ being extended to $U\times U$ by
$$
Y(P(D)u,z)Q(D)v=P(d/dz)Q(D-d/dz)Y(u,z)v,
$$
$$
Dc=0 \qquad\mbox{and}\qquad D=D\otimes 1 \quad\mbox{on}\quad\C[D]\otimes
S'\, ,
$$
and check whether it is a
vertex Lie algebra, i.e. check whether the elements of the form (\ref{881})
and
(\ref{882}) are zero for all $u,v,w\in S'$.

In this case for any given $\ell\in\C$ the quotient
$\V_\ell(U)=\V(U)/\langle c-\ell\,\1\rangle$
of the universal enveloping vertex superalgebra $\V(U)$ is a vertex operator
superalgebra with a conformal vector $\omega\in S_2\subset\V_\ell(U)$.
\end{rem}

{\bf Example 4.}\ \ The only possible examples of (J+ss)-injective formulas
$S$
with a conformal
vector $\omega$ such that $S=S^0=\C c\oplus S_2$ are a special case of
Example 3:
$B=S_2$ is a finite dimensional associative commutative algebra with the
identity
$\omega$ and a symmetric associative bilinear form
such that $\langle\omega,\omega\rangle=1$.
For $u,v\in B$ the formula
$$
Y(u,z)v=\frac{D(u\cdot v)}{z}+\frac{2u\cdot v}{z^2}+
\frac{\frac{1}{2}\langle u,v\rangle c}{z^4}
$$
implies in the corresponding Lie algebra $\L(U)$ the commutation relations
$$
[u(n),v(m)]=(n-m)(u\cdot v)(n+m)
+\frac12{n+1 \choose 3}\delta_{n+m,0}\langle u, v\rangle c.
$$
A somewhat different construction of this Lie algebra and the corresponding
vertex operator algebra is given in [La]. See also [DLM] and the references
therein.


\vspace{1cm}
\noindent Department of Mathematics, University of Zagreb, Bijeni\v{c}ka 30,
10000 Zagreb, Croatia

\noindent e-mail address: primc@cromath.math.hr

\end{document}